\documentclass[12pt]{article}

\usepackage{graphicx,eucal,mathrsfs}
\usepackage{amssymb,amsmath,amsxtra}
\usepackage[lined,linesnumbered]{algorithm2e}
\usepackage{wrapfig}
\usepackage{lipsum}
\usepackage{tikz}
\usetikzlibrary{shapes,arrows}
\usepackage{placeins}
\usepackage{subcaption}

\newcommand{\R}{\mathbb{R}}
\newcommand{\N}{\mathbb{N}}

\newcommand{\E}{\mathbb{E}}

\newcommand{\tr}{\,\textrm{tr}\,}

\begin{document}


\title{{ \Large \sc Shape Gradients for the Failure Probability of a Mechanic Component under Cyclic Loading\\
--- A Discrete Adjoint Approach --- } }
\author{{\sc Hanno Gottschalk} and {\sc Mohamed Saadi}\\
{\small Fachgruppe f\"ur Mathematik und Informatik, Bergische Universit\"at Wuppertal, Germany}\\{\tt \small $\{$hanno.gottschalk,saadi$\}$@uni-wuppertal.de}}

\maketitle

\abstract{\noindent This work provides a numerical calculation of shape gradients of failure probabilities for mechanical components using a first discretize, then adjoint approach. While deterministic life prediction models for failure mechanisms are not (shape) differentiable, this changes in the case of probabilistic life prediction. The probabilistic, or reliability based, approach thus opens the way for efficient adjoint methods in the design for mechanical integrity. In this work we propose, implement and validate a method for the numerical calculation of the shape gradients of failure probabilities for the failure mechanism low cycle fatigue (LCF), which applies to polycrystalline metal. Numerical examples range from a bended rod to a complex geometry from a turbo charger in 3D.}

\vspace{.3cm}

\noindent {\bf Key words:}  Shape Gradients, Failure Probabilities, Adjoint Equation for Structural Mechanics\\
\noindent {\bf Mathematics Subject Classification (2010):} 49Q12, 74P10, 65C50

\section{Introduction}

\begin{figure}[t]
\centerline{
\includegraphics[width=.3\textwidth]{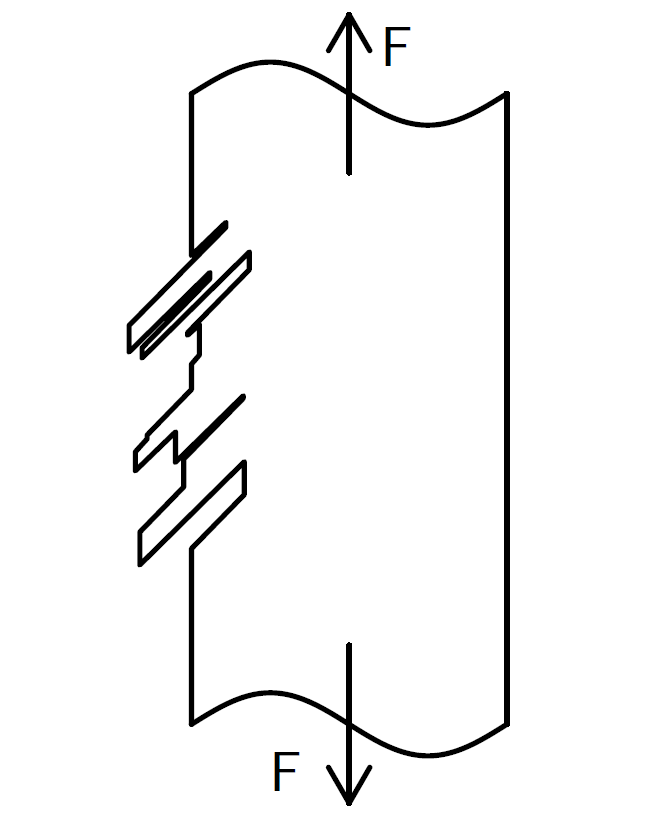}
\includegraphics[width=.5\textwidth]{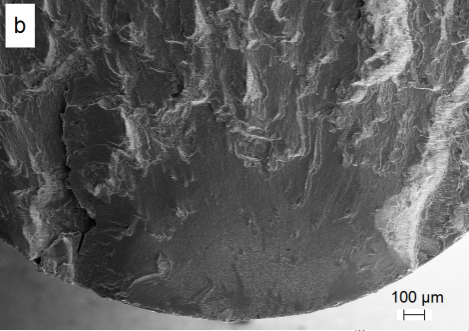}}
\caption{\label{fig:Fatigue} (a) Intrusions and extrusions form on the surface,  (b) a crack initiation from the surface (taken from \cite{GSSKRB}).}
\end{figure}

Mechanical components fail under cyclic loading with stresses well below the ultimate tensile strength (UTS) of the material. This material degeneration process is known as fatigue \cite{BHR,RV}. In polycrystalline metal, the number of identical load cycles until the initiation of a fatigue crack exposes a statistical scatter which, even under controlled lab conditions, is almost an order of magnitude bigger than, say, the 50\% quantile of the data \cite{RV}. The physical origin of the degradation process is the gliding of linear lattice defects (displacements) along crystallographic planes of densest packing in the direction of so-called slip systems \cite{Got}. Under the cyclic reversal of loads, the displacements form intrusions and extrusions at the surface of the component, from which a crack can start to grow, see \cite{Got,BHR,RV} and Fig. \ref{fig:Fatigue}. This is known as low cycle fatigue (LCF) which, in particular, is a surface driven failure mechanism. The random nature of LCF crack formation is then attributed to the random grain structure in a polycrystalline metal. Here, among other effects, the random relative orientation of slip systems relative to the principal stresses and crack percolation over grain boundaries play an important r\^ole \cite{SMB}.

One of the most frequent design objectives in mechanical engineering is mechanical integrity, which here is understood as the absence of cracks. While the loads are calculated with a finite element (FEA) discretization of the (linear) elasticity partial differential equation (PDE), the save life (or design life), i.e. the number of cycles which the component may be safely used, is then calculated for each point (finite element node) at the surface of the component. Here, the local stress values from the FEA are put into some stress-life or strain-life curves, see e.g. \cite{BHR,RV} or Section \ref{sec:FP}. Taking the minimum of all these values -- in combination with safety factors which account for the scatter band in the material curves -- then results in the safe life of a component. The design objective 'mechanical integrity' can then be understood as the task to maximize the safe life.

Automated optimization of mechanical components has become a standard process in mechanical engineering. Here the wish for an efficient calculation of the sensitivity of design objectives in terms of small changes in the geometry has lead to the implementation of (discrete) adjoint equations in commercially available or open source computational fluid dynamics (CFD) solvers. Besides the usefulness of gradients for optimization purposes, the local sensitivities themselves provide also valuable insight to the designers. Unfortunately, the situation up to now is different for the design objective of 'mechanical integrity'.  The reason is properly the deterministic lifing approach, which is based on a (non differential) minimum over a large number of nodes, instead of an integral quantity in the flow variables, as in CFD. This is far more than just a mathematical concern. When this is ignored, as it sometimes happens in the engineering practice, gradient based optimizers tend to jump from one branch of the $\min$-function to the other, if the point of highest loading has a cross-over from one part of the component to some (remote) part of it, so that the convergence process is disrupted.

It has been observed in \cite{GS,BGS,Schmi} that the above three problem statements are interconnected. A proper modeling of the random nature of the LCF crack-initiation process will automatically result in an integral expression in the local stress fields for the acute risk exposure (hazard rate) of the component, which contains the same information as failure probabilities. An adequate model has been proposed in \cite{SRGK1,SRGK2,SSBRKG,LSGB}, see also \cite{HV} and \cite{BSSST} for a related approaches.  Therefore, the passage to a probabilistic design prospective with respect to mechanical integrity, as a by-product, solves the problem of differentiability of objective functionals.   Analytic calculations using the continuous calculus of shapes can be found in \cite{HM,SZ}, see also \cite{ABFJ}.  In this work, we report on the first numerical calculation of shape gradients of the aforementioned probabilistic objective functionals, encoding risk exposure for LCF failure. The method of choice is the 'first discretize, then adjoin' approach which can be realized using standard commercial or open source FEA solvers for the adjoint equation \cite{BR,BS,Troe}. This would not be true for a 'first adjoin, then discretize' approach where derivatives of the original solution up to second order are required, see \cite{SZ}. These are usually not available in the $H^1$-element classes of the commercial solvers. The 'first discretize, then adjoin' approach chosen here therefore seems to us more compatible with industrial application.

Let us underline that it is only the combination of well known ingredients that has lead to this work. Local functionals for the prediction of the UTS of ceramic materials have a history that dates back to Weibull \cite{Wei,BC}. In the design for ceramic materials, these methods are standard in the science of ceramic materials \cite{FM}. Here also the usage of these models in FEA-post-processors with application, e.g., to ceramic heat shields can be found, see e.g. \cite{RBZ}. See also \cite{LLMMW} and articles related to the DARWIN tool suite\footnote{www.darwin.swri.org/} for further probabilistic post processing approaches.

The gradient based, non parametric optimization of shapes, including topology optimization, is a well established and advanced mathematical field, both from the analytical, see \cite{SZ,HM,Troe}, and the numerical side, see \cite{All,BS,BR} for incomplete lists. The impressive results obtained here in the case of mechanical design are mostly linked to another design objective, namely the global stiffness of the component, encoded in the compliance as objective functional. This stands for a design intent which is different from mechanical integrity.

Our paper is organized as follows: In Section \ref{sec:FP} we review the classical LCF life prediction and the related probabilistic model. Section \ref{sec:DIS} explains the implementation of this model as a FEA post processor. Here we essentially follow \cite{SRGK1,SRGK2,Schmi}. Section \ref{sec:AE} contains the derivation of the shape gradients on the basis of the discrete adjoint equation. Some lengthy calculations have been transferred to Appendix \ref{sec:AppCalc}. Section \ref{sec:NUM} discusses some detailed information on the implementation and presents two numerical test cases: firstly,  a bended rod exposed to a tensile load and secondly a complex turbo geometry (taken from \cite{Wit}). Finally, in Section \ref{sec:OUT}, we draw our conclusions and give an outlook to future work. Some of our results have been announced in \cite{GSDKS}.

\section{Failure Probabilities for Low Cycle Fatigue}
\label{sec:FP}

\subsection{The Elasticity PDE}
Let $\Omega\subseteq \R^3$ a bounded region filled with some polycrystalline metal. This represents the shape of the component initially in force free equilibrium. By $H^1_D(\Omega,\R^3)$ we denote the Sobolev space of $L^2(\Omega,\R^2)$ functions $u$ with weak 1st derivatives $\nabla u\in L^2(\Omega,\R^{3\times 3})$ such that $u=0$ on $\partial\Omega_D$ \cite{Cia}. Here $\partial\Omega$ is a piecewise Lipschitz boundary such that the part with Dirichlet boundary $\partial\Omega_D$ is an open portion of $\partial\Omega$ with non vanishing surface measure. Furthermore we set the part of natural boundary conditions to be $\partial\Omega_N=(\partial\Omega\setminus\partial\Omega_D)^\circ$, where $^\circ$ stands for the open interior in $\partial\Omega$.

The function $f\in L^2(\Omega,\R^3)$ stands for a volume force density like imposed by gravity or centrifugal loads. $g\in L^2(\partial\Omega_N,\R^3)$ is a surface load, like caused by static gas pressure $P(x)$, $x\in \partial\Omega_N$, $g(x)=-P(x)n(x) $ with $n(x)$ the outward normal vector field on $\partial\Omega_N$. Let $u\in H^1(\Omega,\R^3)$ be the displacement field caused by these loads in a way that will be specified shortly. We define $\varepsilon(u)=\frac{1}{2}(\nabla u+\nabla u^T)\in L^2(\Omega,\R^{3\times 3})$ the linearized strain tensor field and we define the stress
tensor field $\sigma \in L^2(\Omega,\R^{3\times 3})$ via $\sigma(u)=\lambda \tr(\varepsilon(u))\,{\cal I}+2\mu \, \varepsilon(u)$, where $\lambda,\mu>0$ are the Lam\'e coefficients that represent the macroscopically isotropic elastic behavior of the material.  ${\cal I}$ is the unit matrix on $\R^3$, $\tr$ is the trace on the space of $3\times 3$ matrices and $M^T$ stands for the transpose of a matrix $M$.

It is well known  that under the given conditions, the PDE system of linear elasticity in the weak form
\begin{equation}
\label{eqa:WeakPDE}
   B(u,v)=\int\limits_{\Omega}f\cdot v\, dx+\int\limits_{\partial\Omega_N}g\cdot v \,dA , \forall v\in H^ 1_{D}(\Omega,\R^ 3)
\end{equation}
has an unique solution $u\in H^1(\Omega,\R^3)$ with the bilinear form given by
\begin{equation}
B(u,v)=\int\limits_{\Omega}\sigma(u):\varepsilon(v)\, dx=\lambda\int\limits_{\Omega} \nabla\cdot u\nabla\cdot v\,dx+2\mu\int\limits_{\Omega}\varepsilon(u):\varepsilon(v)\, dx.
\end{equation}
Here $\cdot$ stands for the scalar product on $\R^3$ and $\varepsilon(u):\varepsilon(v)=\tr(\varepsilon(u)\varepsilon(v))$ is the contraction of both indices of $\varepsilon(u)$ and $\varepsilon(v)$, $dx$ is Lebesgue measure and $dA$ is the surface measure on $\partial\Omega$. In fact this follows from Korn's \cite{Cia} inequality which for the present case gives the coercivity $A^{-1}\|u\|^2_1\leq B(u,u)\leq A\|u\|^2_1$ $\forall u\in H^1_D(\Omega,\R^3)$ for some constant $\infty>A>0$ which is needed to apply the Lax-Milgram theorem.

In the following, some operations will frequently require higher degrees of regularity of the solution $u$ than merely being a $H^1$-Sobolev function. In particular this applies to the pointwise evaluation of the stress tensor in deterministic LCF life prediction or the restriction of the stress tensor to the boundary in the local, probabilistic model for LCF life. Both operations  can be made rigorous using elliptic regularity theory, see e.g. \cite{AGN,GS}  for suitable conditions on boundary regularity and  of the loads $f$ and $g$.

\subsection{The Equations of Ramberg-Osgood and Coffin-Manson Basquin}
In this section we briefly review the deterministic LCF life prediction procedure which forms the basis of the probabilistic model. We closely follow the materials science literature \cite{BHR,Got,RV}.

\begin{figure}[ht]
\centerline{
\includegraphics[width=.32\textwidth]{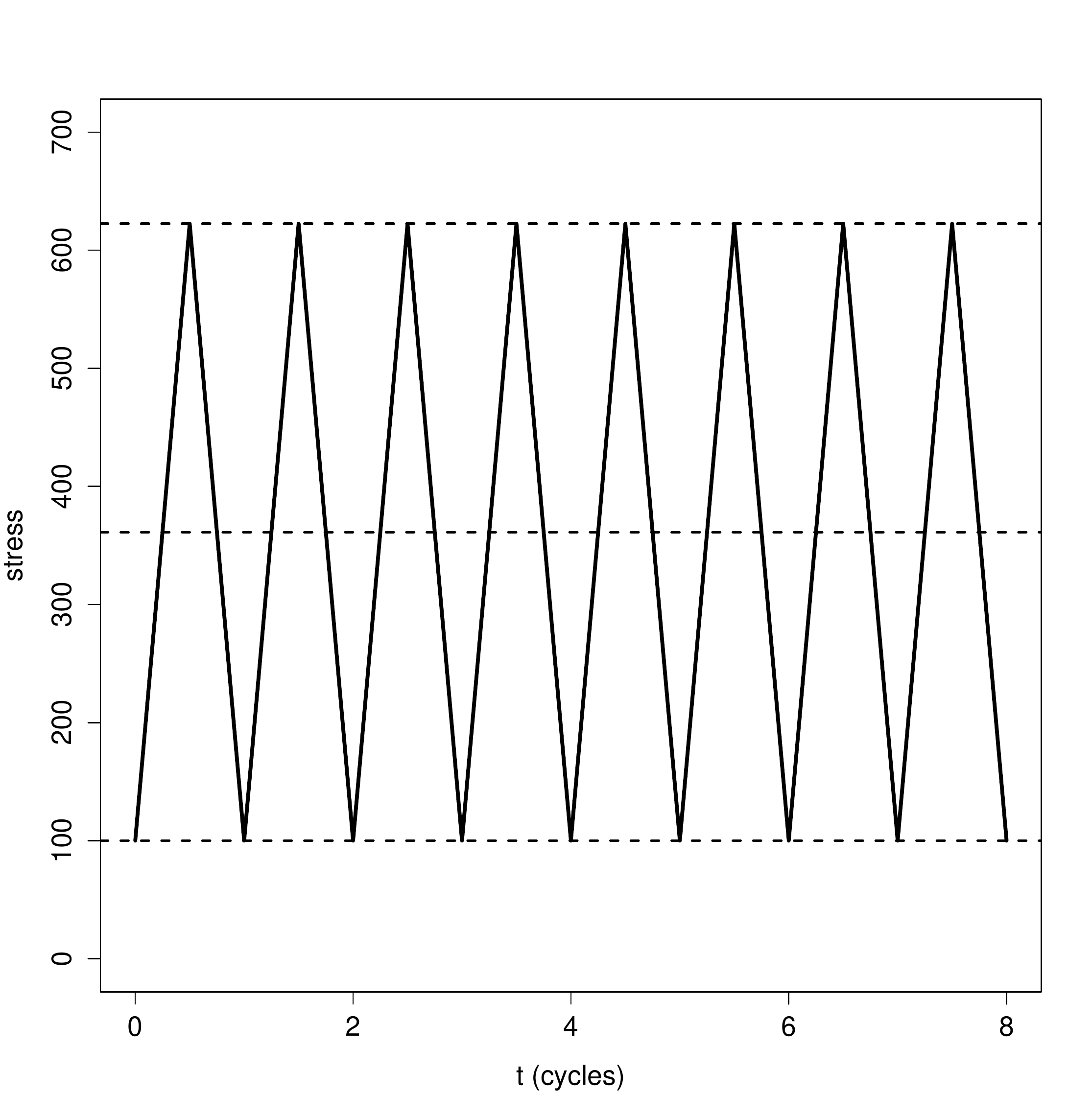}
\includegraphics[width=.32\textwidth]{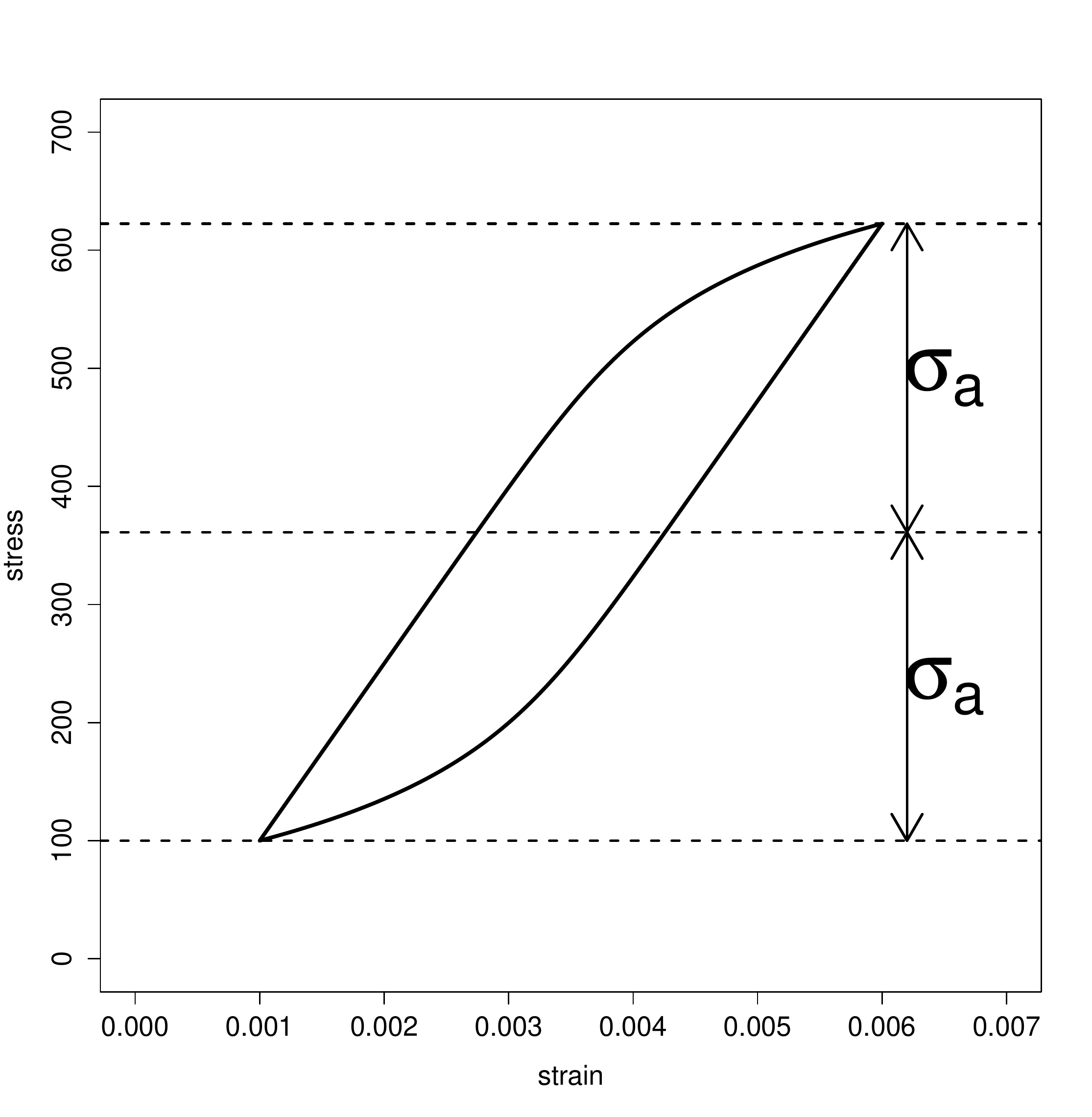}
\includegraphics[width=.32\textwidth]{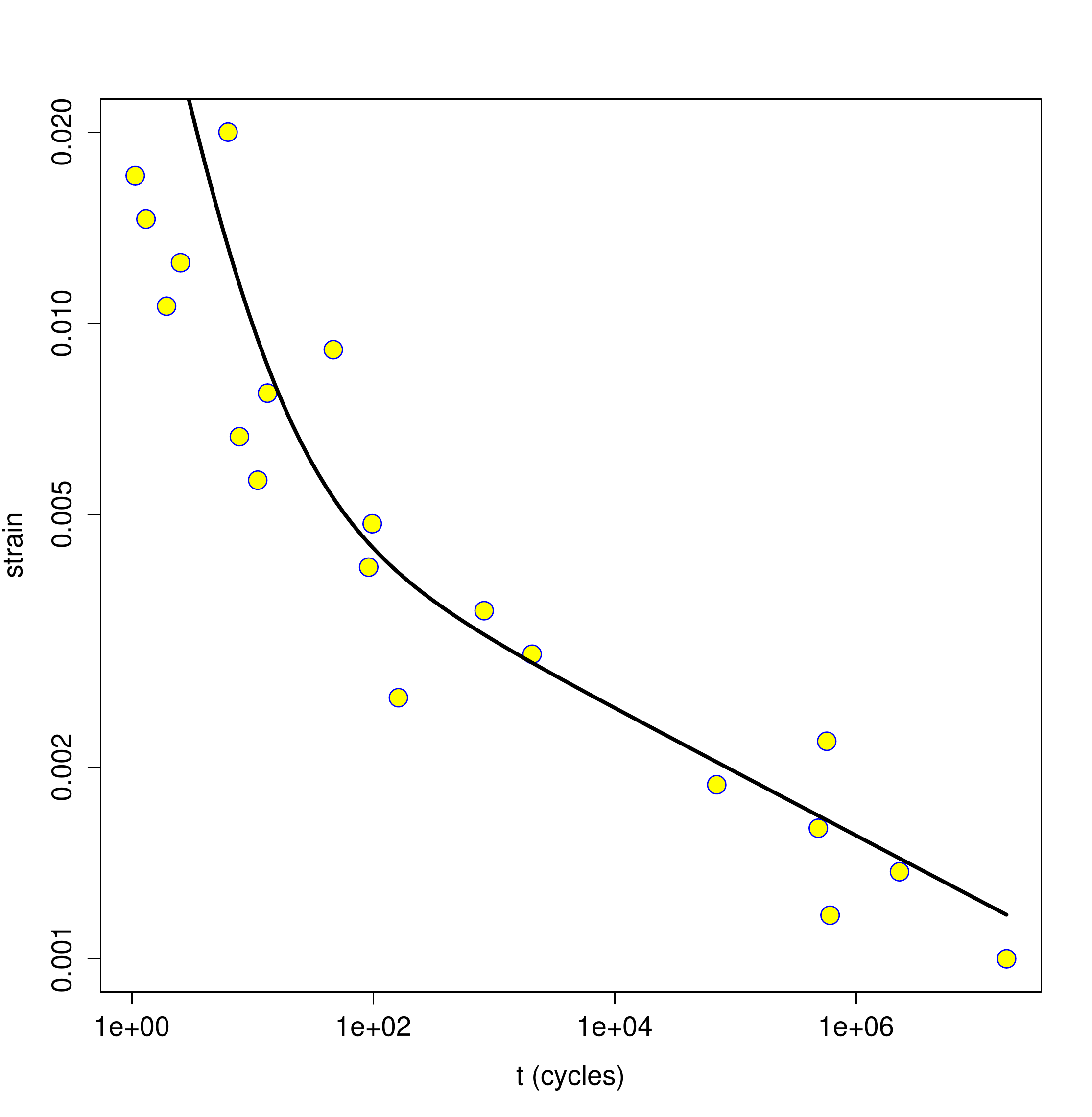}}
\caption{\label{fig:ROCMB} (a) Load cycles, (b) elastic-plastic stress-strain hysteresis with branches according to the Ramberg-Osgood   equation, stress in MPa, (c) Coffin-Manson-Basquin curve fitted to  specimen test results (data simulated on the basis of realistic values). Note that both axes are on $\log$-scale.}
\end{figure}

For $x\in\partial\Omega$ let $\sigma=\sigma(u)$ be the stress tensor from the solution of \eqref{eqa:WeakPDE} evaluated at $x$. For simplicity we assume that the other extreme of the load cycle is the unloaded situation $u=0$. If this is not the case, we can make use from the linearity of the elasticity PDE and set $f=f_1-f_2$ and $g=g_1-g_2$, where $f_j$ and $g_j$, $j=1,2$, are the volume and surface loads at the extremes of the load cycle, see Figure \ref{fig:ROCMB} (a). As hydrostatic stress can not induce shear stress\footnote{The slip directions $s$ mentioned in the introduction are perpendicular to the normal $\nu$ of the plane of densest packing, thus $s\cdot \sigma\nu=0$ if $\sigma=\sigma_0{\cal I}$, $\sigma_0\in\R$, and no force acts on the slip system $s$ \cite{Got}.}, we project $\sigma$ to the trace-free matrices via $\sigma'=\sigma-\frac{1}{3}\tr (\sigma){\cal I}$.  In the next step we use the Hilbert Schmid norm $\sigma_v^2=\frac{3}{2}\sigma':\sigma'\in\R_+$ to define the elastic von Mises comparison stress. The normalization is chosen in order to reproduce the first principal stress in the case of an uniaxial stress state.  Setting $\sigma_a=\frac{1}{2}\sigma_v\in\R_+$ we arrive at the so called elastic amplitude stress.

Let us briefly consider elastic-plastic material behavior based on a load cycle with an elastic-plastic amplitude stress $\sigma_a^{\rm el-pl}$. We postpone a mathematical definition for a short moment. Then, the elastic-plastic stress amplitude is connected to the elastic-plastic strain amplitude via the Ramberg-Osgood (RO) relation \cite{RO}
\begin{equation}
\label{eqa:RO}
\varepsilon_a^{\rm el-pl}={\rm RO}(\sigma^{\rm el-pl}_a)=\frac{\sigma_a^{\rm el-pl}}{E}+\left(\frac{\sigma_a^{\rm el-pl}}{K}\right)^{1/n'},
\end{equation}
see also Fig. \ref{fig:ROCMB} (b). $E=\frac{\mu(3\lambda+2\mu)}{\lambda+\mu}$ is Young's modulus, $K$ is the hardening coefficient and $n'$ is the hardening exponent.

Finally we use the Coffin-Manson Basquin (CMB) equation to predict the component's LCF-life, i.e. the number of cycles to crack initiation ${\rm Ni}_{\rm det}$, via
\begin{equation}
\label{eqa:CMB}
\varepsilon^ {\rm el-pl}_a={\rm CMB}({\rm Ni}_{\rm det})=\frac{\sigma_f'}{E}\left(2{\rm Ni}_{\rm det}\right)^b+\varepsilon_f'\left(2{\rm Ni}_{\rm det}\right)^c.
\end{equation}
 Here, $\varepsilon_f',\sigma_f'>0$ and $b<c<0$ are material constants that are determined from experimental data, see Fig.\ \ref{fig:ROCMB} (c).

Let us now return to the definition of the elastic-plastic stress amplitude $\sigma_a^{\rm el-pl}$ on the basis of the stress tensor field $\sigma=\sigma(u)$  obtained from the solution of \eqref{eqa:WeakPDE}. In the engineering practice, this is often done with elastic-plastic stress conversion rules by Neuber and Glinka \cite{Neu,HS}.
Let us, for simplicity, explain the Neuber shake-down rule, only. It is obtained by equating the (uniaxial) elastic-plastic energy density $\varepsilon_a^{\rm el-pl}\sigma_a^{\rm el-pl}$ and the corresponding purely elastic quantities $\sigma_a\varepsilon_a=\frac{\sigma^2_a}{E}$. As $\sigma_a$ is calculated on the basis of the elasticity PDE, we can now use the Ramberg-Osgood equation to obtain
\begin{equation}
\label{eqa:ShakeDown}
\sigma_a={\rm SD}(\sigma_a^{\rm el-pl})=\sqrt{E\left(\frac{(\sigma_a^{\rm el-pl})^2}{E}+\sigma_a^{\rm el-pl}\left(\frac{\sigma_a^{\rm el-pl}}{K}\right)^{1/n'}\right)}.
\end{equation}
This gives the following formula for the deterministic life prediction:
 \begin{equation}
 \label{eqa:Ndet}
 {\rm Ni}_{\rm det}(\sigma(x))={\rm CMB}^{-1}\circ \rm{RO}\circ{\rm SD}^ {-1}(\sigma_a(x))),~~{\rm Ni}_{\rm det}=\inf_{x\in\partial \Omega}{\rm Ni}_{\rm det}(\sigma(x)).
 \end{equation}
We note that there exist various variants of the above general scheme - e.g. one might prefer to apply the shake down to the von Mises stress $\sigma_v$ instead of the amplitude stress.  The numerical differences, especially after fitting the parameters to identical experimental data sets, in most cases are  marginal. In real design applications, safety factors are used in various parts of the procedure in order to account for the scatter in the experimental data and modeling errors.

\subsection{The Local Probabilistic Model for LCF Failure}

In the local, probabilistic model for LCF \cite{GS,Schmi,SRGK1,SRGK2,SSBRKG} crack formation is modeled as a random event in space, i.e. in $x\in \partial\Omega$, and time $\geq 0$ counting the number load cycles that the component has undergone until the formation of a crack. Although $t$ strictly speaking is a integer number, already in \eqref{eqa:CMB} we tacitly passed to continuous time.

Let thus $[s,t]$, $0\leq s<t$, be a time interval and let $A\subseteq \partial\Omega$ be some surface region. We then denote by $\gamma([s,t],A)$ the random number of cracks that initiated in $A$ during the time interval $[s,t]$. In mathematical terms, $\gamma$ is a so-called point process \cite{Kal}. If we assume that cracks in non-intersecting regions $A_1,\ldots,A_n$ of $\partial\Omega$ this initial phase do not interact with each other, i.e. $\gamma([0,t],A_1)$,\ldots,$\gamma([0,1],A_n)$ are independent random variables, we can deduce under quite general assumptions that $\gamma([0,t],A)$ must be Poisson point process (PPP) \cite{Kal,Wat}. In particular, there exists a sigma finite measure $\rho(t,.)$ on $\partial\Omega$ such that
\begin{equation}
\label{eqa:Pois}
P(\gamma([0,t],A)=n)=e^{-\rho(t,A)}\frac{\rho(t,A)^n}{n!},~~n\in\N.
\end{equation}
Here $P$ is the underlying probability measure. We also note that $\rho(t,A)$ is the expected number of cracks that originated in $A$ up to time $t$, i.e. $\E[\gamma([0,t],A)]=\rho(t,A)$.

Let us shortly review some modeling principles that lead to the model. First of all, the expected number of cracks $\rho(t,A)$ should only depend on the local stress state $\sigma(x)$ on $A$. Furthermore $\rho(t,A)$ has to be monotonically increasing in $t$ as cracks, one they initiated, will stay in place. This leads to
\begin{equation}
\label{eqa:LocalForm}
\rho(t,A)=\int\limits_0^t\int\limits_{A}\varrho(t,\sigma(x))\, dA\, dt\,,
\end{equation}
where we omitted the possibility that the local crack formation intensity up to time $t$ could also depend on derivatives of the stress tensor field.

We now need a reasonable model for the local crack formation intensity $\rho(t,\sigma(x))$. On the one hand, this model should reflect the experimentally well established procedure of deterministic LCF-life prediction as a kind of average (or quantile) value. On the other hand, the model should allow to adjust itself to different levels of statistical scatter in LCF-life times. We choose the following
\begin{equation}
\label{eqa:locWeibull}
\varrho(t,\sigma(x))=\frac{\bar{m}}{{\rm Ni}_{\rm det}(\sigma(x))}\left(\frac{t}{{\rm Ni}_{\rm det}(\sigma(x))}\right)^{\bar m-1},
\end{equation}
where ${\rm Ni}_{\rm det}(\sigma(x))$ as defined in \eqref{eqa:Ndet} is the scale parameter and $\bar{m}>0$ is the shape parameter in the sense of the scale-shape class of probability distributions, see e.g. \cite{EM} for a discussion with application to reliability statistics. In the limit $\bar m\to\infty$, we recover the deterministic model. We then get,
\begin{equation}
\label{eqa:MeanFailure}
\rho(t,A)=\int\limits_{A}\left(\frac{t}{{\rm Ni}_{\rm det}(\sigma(x))}\right)^{\bar m} dA \, .
\end{equation}
We define the surface integral type objective functional as
\begin{equation}
\label{eqa:Obj}
J(\Omega,u)=\int\limits_{\partial\Omega}\left(\frac{1}{{\rm Ni}_{\rm det}(\sigma(x))}\right)^{\bar m} dA\, .
\end{equation}
 We now deduce the probability of failure ${\rm PoF}(t)$ -- understood as the presence of a LCF crack -- as a function of time (load cycle count) $t$ 
\begin{equation}
\label{eqa:PoF}
{\rm PoF}(t)=1-S(t)=1-P(\gamma([0,t],\partial\Omega)=0)=1-e^{-t^{\bar{m}}J(\Omega,u)},~t\geq 0.
\end{equation}
In \eqref{eqa:PoF}, $S(t)$ is the probability of survival, which is defined here as the absence of cracks up to time $t$ on the entire component's surface $\partial\Omega$, hence the event $\gamma([0,t],\partial\Omega)=0$. Equation \eqref{eqa:PoF} now follows from \eqref{eqa:Pois}, \eqref{eqa:MeanFailure} and \eqref{eqa:Obj}. We recognize \eqref{eqa:PoF} as a Weibull distribution with scale parameter $\left(\frac{1}{J(\Omega,u)}\right)^{\frac{1}{\bar m}}$ and shape parameter $\bar{m}$.

Note however that the CMB constants $\sigma_f'$ and $\varepsilon_f'$ are not the same as in the deterministic life calculation, see \cite{SSBRKG} for the details on the calibration of the probabilistic model, see also Subsection \ref{sec:DISEx} below.  We also note that $\bar m>1$ is needed for situation when the hazard rate, which is roughly the probability rate that failure will happen in the next moment, given the component survived until that moment, grows with the time $t$. See \cite{EM} for further details.

\section{Discretization with Finite Elements}
\label{sec:DIS}

\subsection{Discretization of the Probabilistic Model}

 The numerical approximation of the failure probability ${\rm PoF}(t)$ in \eqref{eqa:PoF} requires the discretization of the objective functional \eqref{eqa:Obj}, which is given as a surface integral. Let $\mathcal{T}_h$ be a set of finite elements with reference element $\widehat{K}$ and let $T_K:\widehat{K}\to K$ be the standard transformation map from the reference element to $K\in \mathcal{T}_h$. See Appendix \ref{app:FE} and \ref{app:DiscretePDE} for further details. Let thus $\mathcal{N}_h$ be the collection of all faces $F$ of finite elements
 $K=K(F)\in\mathcal{T}_h$ that lie in $\partial\Omega$. The total number of such faces is denoted by $N_{\rm fc}^F$. Let $\widehat{F}$ be the face in the reference element $\widehat{K}$ such that $T_K:\widehat{F}\to F$.
  We chose quadrature points $\hat \xi_l^F$ and weights $\hat\omega_{l}^F$ on $\widehat{F}$. We then get from \eqref{eqa:Obj}, if we insert the approximate solution $u$ of \eqref{eqa:WeakPDEDiscrete},
\begin{align}
 \label{eqa:ObjDiscretized}
\begin{split}
 J(\Omega,u)&=\sum_{F\in\mathcal{N}_h}\int\limits_F \left(\frac{1}{{\rm Ni}_{\rm det}(\sigma(x))}\right)^{\bar m} \,dA\\
 &=\sum_{F\in\mathcal{N}_h}\int\limits_{\widehat{F}} \left(\frac{1}{{\rm Ni}_{\rm det}(\sigma( T_{K(F)}(\hat x)))}\right)^{\bar m} \sqrt{\det g_F(\hat x)} \,d\hat A \\
 &\approx \sum_{F\in\mathcal{N}_h}\sum_{l=1}^ {l_q^ F}\hat \omega_l^F\left(\frac{1}{{\rm Ni}_{\rm det}(\sigma(T_{K(F)}(\hat \xi_l^ F)))}\right)^{\bar m} \sqrt{\det g_F(\hat\xi_l^ F)}.
 \end{split}
\end{align}
Here $\det g_F(\hat\xi)$ is the Gram determinant, i.e. the determinant of
\begin{equation}
\label{eqa:Gram}
 g_F(\hat\xi)= \widehat{\nabla}^F(T_K\restriction \widehat{F})(\xi)\left(\widehat{\nabla}^F(T_K\restriction \widehat{F})\right)^T(\hat \xi),~~\hat\xi\in\widehat{F	},
 \end{equation}
with $\widehat{\nabla}^F$ the (component wise) gradient on $\widehat{F}$. Here $\sigma(x)=\lambda \nabla\cdot u(x)\mathcal{I}+2\mu \varepsilon(u(x))$ can be found by expanding $u$  in the global basis functions $\theta_j$, we get for $u(\xi)$, 
$u(\xi)=\sum_{j=1}^{N}u_j\theta_j(\xi)=\sum_{K\in\mathcal{T}_h}\sum_{m=1}^{n_{\rm sh}}u_{j_{(K,m)}}\widehat{\theta}_m \circ T_{K}^{-1}(\xi)$ so that
\begin{equation}
\label{eqa:Jacobi}
\nabla u(\xi)=\sum_{m=1}^{n_{\rm sh}}u_{\widehat{j}(K,m)}\otimes(\widehat{\nabla} T_{K}(\widehat{\xi})^{T})^{-1}\widehat{\nabla} \widehat{\theta}_m (\widehat{\xi})),~\mbox{ for } \xi \in K\mbox{ and }\xi=T_K(\hat\xi).
\end{equation}
Symmetrizing the above Jacobi matrix, we obtain $\varepsilon(u(\xi))$.  Finally, we have
\begin{equation}
\label{eqa:Divergence}
\nabla\cdot u(\xi)=\sum_{m=1}^{n_{\rm sh}}\tr\left(u_{\widehat{j}(K,m)}\otimes(\widehat{\nabla}T_{K}(\widehat{\xi})^{T})^{-1} \widehat{\nabla}
\widehat{\theta}_m(\widehat{\xi}))\right) \mbox{ for } \xi \in K\mbox{ and }\xi=T_K(\hat\xi).
\end{equation}
Note that $\widehat{\nabla} T_{K}(\hat\xi)=\sum_{j=1}^{n_{\rm sh}}\widehat\nabla\widehat\theta_j(\widehat{\xi})X_{K,j}$ is easily calculated from \eqref{eqa:Transformation}. For the discretized objective functional, we also use the notation $J(X,U)$, where $X$ and $U$ the $N\times 3$-tensor of global coordinates and global degrees of freedom, respectively.

\subsection{Numerical Validation and a First Application}
\label{sec:DISEx}
\begin{figure}[ht]
\centerline{
\includegraphics[width=.43\textwidth]{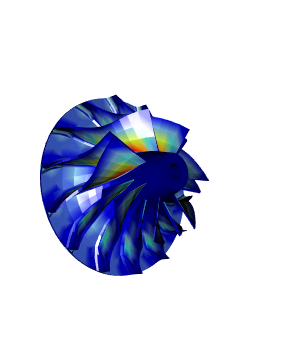}
\includegraphics[width=.53\textwidth]{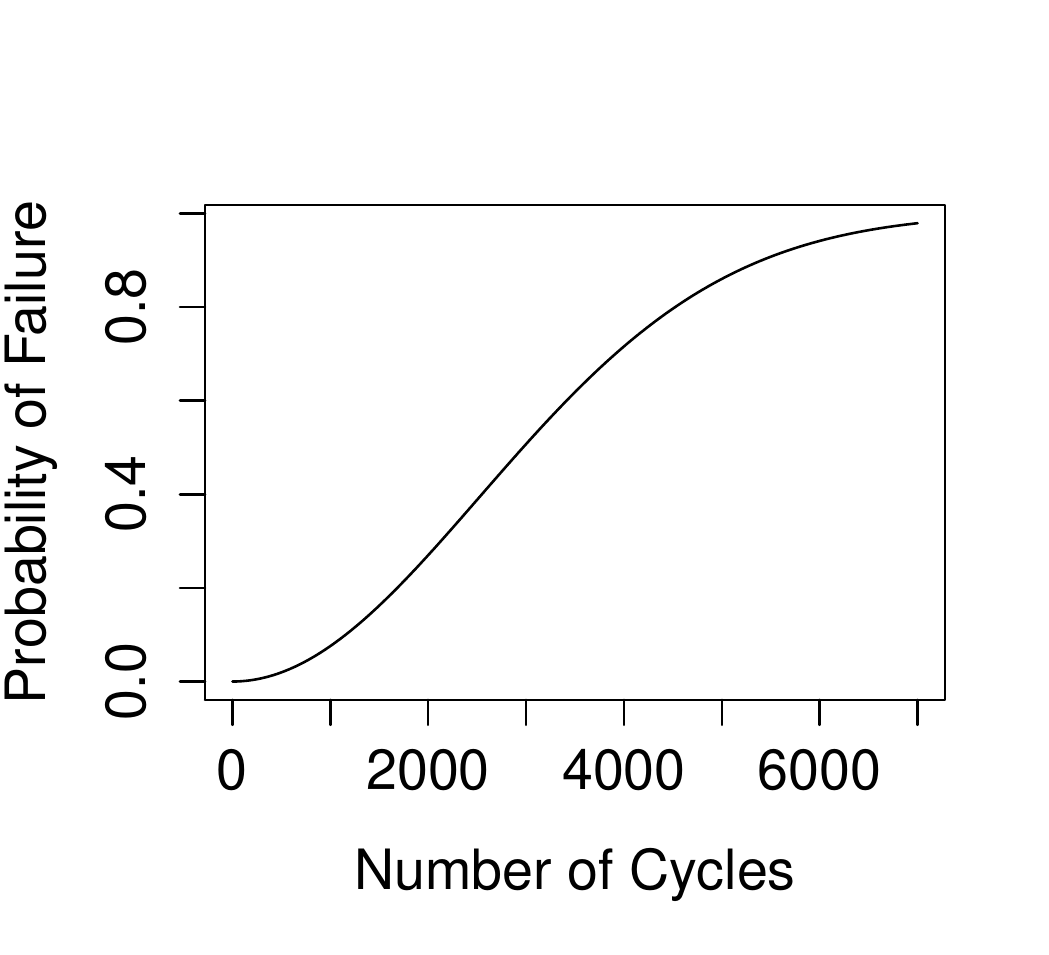}
}
\caption{\label{fig:Compressor} Probability of failure over the number of load cycles (left) for a radial compressor model. (a) crack initiation intensity $\rho(1,\sigma(x))$ in logarithmic color code from blue (low) to red (high). (b) PoF over the number of load cycles for the compressor model. }
\end{figure}

As a first application, we present the calculation of crack initiation PoF over the number of load cycles for a radial compressor of a model gas turbine, see \cite{Wit} for the documentation. The model is part of the CalculiX open source FEA tool suite\footnote{www.calculix.de}, which is a design tool for MTU Aero Engines. The outer diameter is $87$ mm and the construction material is taken to be the alloy AlMgSi6082, containing roughly $97$ wt\% Al, roughly $1$ wt\% Si and Mg along with minor contributions of Mn, Fe and Cr. The  detailed chemical composition is described in \cite{BACF}.  

Air pressure is neglected, hence only volume forces resulting from the centrifugal load from a rotation speed of  $110\,000$ rpm with a specific density of the material of $2.65$ g/cm${}^3$. The Lam\'e coefficients for the material are $\lambda=\frac{E\nu}{(1+\nu)(1-2\nu)}=40\,385$ MPa and $\mu=\frac{E}{2(1+\nu)}=26\,923$ MPa are calculated from Young's modulus ($E=70\, 000$ MPa) and Poisson's ratio $\nu=0.3$. 

The compressor consists of 7 symmetric segments with 2 slightly different blades each, so that the actual FEA model only represents one seventh of the construction with cyclic boundary conditions. Dirichlet boundary conditions are set at the bore. The FEA mesh consists out of $47\,971$ nodes with $1\,302$ brick elements of type C3D20R with 20 degrees of freedom and reduced quadrature $l_q=8$. The surface quadrature points are selected such that $l_q^F=16$, confer \cite{SRGK1} for a numerical convergence study in the number of quadrature points for a related application.

Ramberg-Osgood Parameters $n'=0.064$ and $K=443.9$ are reported in \cite{BACF}. The CMB parameters are given in the top row of Table \ref{tab:CMB}. The passage from deterministic CMB parameters \cite{BACF} to probabilistic the CMB parameters needed here is performed in two steps. 

We interpret the reported deterministic values as belonging to a $50\%$ quantile curve. This is justified, as in engineering these curves are mostly obtained as least square fits to the logarithmic specimen data in both dimensions $\varepsilon_a^{\rm le-pl}$ and ${\rm Ni}_{\rm det}$. This fit corresponds to a maximum likelihood fit with log-Normal residuals, where the average curve for the logarithmic data corresponds to a 50\%-quantile curve in the non logarithmic data. Within our Weibull model, the relation between the 50\% quantile $q_{0.5}$ and the scale variable $\eta$ is $\eta=q_{0.5}/\log(2)^\frac{1}{\bar{m}}$, and we use this to rescale the CMB parameters according with a factor $\log(2)^{-\frac{b}{\bar{m}}}$ for $\sigma_f'$ and $\log(2)^{-\frac{c}{\bar{m}}}$ for $\varepsilon'_f$, see \cite{SSBRKG}.  The result is given in the second row of  Table \ref{tab:CMB}. 

In the next step, we account for statistical size effect, i.e. the effect that larger structures have the tendency to fail earlier than smaller ones.
In order to convert the CMB parameters that correspond to the scale variable curve of the specimen to the probabilistic CMB parameters of a fictitious specimen with unit surface, we have to apply the scaling relation 
\begin{align}\label{eqa:interpretation_CMB.7.1}
\frac{\sigma_f'(|\partial\Omega_1|)}{\sigma_f'(|\partial\Omega_2|)}=\left(\frac{|\partial\Omega_1|}{|\partial\Omega_2|}\right)^{\frac{b}{\bar{m}}},\quad
\frac{\varepsilon_f'(|\partial\Omega_1|)}{\varepsilon_f'(|\partial\Omega_2|)}=\left(\frac{|\partial\Omega_1|}{|\partial\Omega_2|}\right)^{\frac{c}{\bar{m}}},
\end{align}
see \cite{SSBRKG,Schmi,RV}. The result is given in the third row of Table \ref{tab:CMB}. The value for $|\partial\Omega|=1$ corresponds to the probabilistic parameter in \eqref{eqa:CMB}. Note that due to \eqref{eqa:interpretation_CMB.7.1} this depends on the length scale, which is [mm] in this work (in agreement with CAD models in turbo machinery applications).   The value for $\bar m$ for AlMgSi6082 is not reported in the literature. However, experience with steel and Ni-based superalloys shows that $\bar m$ is in the region of $2$ for a large range of materials.  

\begin{table}[t]
\centerline{
\begin{tabular}{||l|c|c|c|c|c|c||}
\hline
&$\sigma'_f$(MPa)&$\epsilon_f'$&$b$&$c$&$\bar m$&$|\partial\Omega|$ (mm${}^ 2$)\\\hline
Deterministic &$487$&$0.209$&$-0.593$&$-0.07$&--&$377$\\
($q_{0.5}$-curve)&&&&&&\\\hline
Weibull scale  &$436$&$0.206$&$-0.593$&$-0.07$&$2$&$377$\\
($q_{0.63}$-curve)&&&&&&\\\hline
Probabilistic &$2536$ &$0.254$ &$-0.593$&$-0.07$&$2$&$1$\\
&[MPa$\times$ mm${}^{\frac{2b}{\bar m}}$]&[ mm${}^\frac{2c}{\bar m}]$&&&&\\
\hline
\end{tabular}}
\caption{\label{tab:CMB} Material parameters for the probabilistic model and their deterministic counterparts (taken from \cite{BACF}). The parameter $\bar{m}$ is a guess on the basis of probabilistic investigation of other polycrytalline metals.   }
\end{table}
   
The FEA model is solved for $u$ with the commercial solver ABAQUS 6.1 on a laptop with Intel Core i7-3632QM CPU @ 2.20 GHz with 12GB RAM with a time consumption of 28 seconds. Peak v. Mises stress at the bore is roughly $310$ MPa. The surface quadrature for the calculation of $J(\Omega,u)$ is calculated with an script in the R 3.2.1 and consumed 19 sec execution time on a single core. The obtained $J$-value has to be multiplied by 7 to account for the seven rotated sectors.

The resulting total $J$- value is $7.8541 \times 10^{-8}$ [cycles${}^{-\bar m}$], which corresponds to a Weibull scale variable $\eta=3568$ cycles. For a Weibull distribution, this corresponds to the $1-\frac{1}{e}\approx 63\%$ quantile. Figure \ref{fig:Compressor} (a) visualizes the local failure intensity \eqref{eqa:locWeibull} and (b) gives the integrated probability of failure (PoF) over the number of cycles $t$. The visual inspection shows that not only the bore region, but also the fillet areas contribute significantly to the total value of $J$. The safe life of $2\,000$ cycles reported in the CalculiX tool suite corresponds to a PoF of roughly $25\%$, according to the local, probabilistic model for LCF, which nicely fits into the over all picture for given design example.

\section{The Shape Sensitivity for the Probability of Failure via the Adjoint Equation}
\label{sec:AE}

\subsection{Straight Forward Calculation of the Shape Sensitivity}
We first fix some conventions. The mesh nodes $X=\{X_1,\ldots,X_N\}$ represent the discretized geometry $\Omega$. Furthermore we recall that $U=(u_j)_{j=1,\ldots,N}$ stands for the global degrees of freedom and $u=\sum_{j=1}^Nu_j\theta_j(x)$ with $\theta_j$ the global shape functions. In this section we equivalently write $J(X,U)$ instead of $J(\Omega,u)$. Note that $U$ is a $N\times 3$ tensor as well, provided we consider the constraints $u_j=0$ for $X_j\in\partial\Omega_D$ as a part of the discretized equation \eqref{eqa:WeakPDEDiscrete},  see also Appendix \ref{app:DiscretePDE}.

$X$ is a $N\times 3$ tensor and, by definition, for a tensor of scalar $Q(X,Z)$, the partial derivative $\frac{\partial Q(X,U)}{\partial X}$ is an augmented tensor with two additional slots of dimension $N$ and $3$. Here $U$ stands for some other tensor with unspecified dimensions (or a scalar, if there are none). If $U$ is an expression  that depends on $X$, $U=U(X)$,  the total derivative of $Q(X,U(X))$ with respect to $X$ is $\frac{dQ(X,U(X))}{dX}=\frac{\partial Q(X,U(X))}{\partial X}+\frac{\partial Q(X,U(X))}{\partial  U}\frac{\partial U(X)}{\partial X}$. Here the last expression has to be understood in the sense that the additional tensor slots generated by the partial $U$-differentiation are contracted with the $U$-slots in $\frac{\partial U(X)}{\partial X}$.

It is now straight forward to see that the probability of failure ${\rm PoF}(t)={\rm PoF}(t,X,U(X))$, see \eqref{eqa:PoF}, at a number of cycles $t\geq 0$ can be differentiated with respect to the geometry $X$ as follows
\begin{equation}
\label{eqa:dPOF_dX}
\frac{d\, {\rm PoF}(t,X,U(X))}{dX}=t^{\bar{m}}\frac{d J(X,U(X))}{dX} \, e^{-t^{\bar{m}}J(X,U(X))}.
\end{equation}
The above equation is already remarkable in the sense that the direction of the gradient with respect to the geometry does not depend on the number of load cycles $t>0$, which excludes the use of the formalism in design to life activities. We now continue the differentiation and obtain
\begin{equation}
\frac{d J(X,U(X))}{dX}=\frac{\partial J(X,U(X))}{\partial X}+\frac{\partial J(X,U(X))}{\partial  U}\frac{\partial U(X)}{\partial X},
\end{equation}
where on both sides we have $N\times 3$ tensors according to our conventions. Note that in the last term on the right hand side $U$-degrees of freedom that correspond to Dirichlet boundary conditions do not depend on $X$ and thus give zero in $\frac{\partial U(X)}{\partial X}$.

In the remainder of this subsection we want to discuss the calculation of $\frac{\partial J(X,U)}{\partial X}$ and $\frac{\partial J(X,U)}{\partial U}$. Inspection of \eqref{eqa:ObjDiscretized} shows that it is easier to deal with $\frac{\partial J(X,U)}{\partial U}$ as here only $\sigma$ depends on $U$, so we do this first. The following algorithm describes the necessary steps. Detailed calculations can be found in Appendix \ref{app:dJ_dU}. Note that here it suffices to calculate all partial derivatives with respect to degrees of freedom that belong to $\mathcal{T}_h^F$, which is the set of finite elements $K\in \mathcal{T}_h$ with at least one surface face $F\in \mathcal{N}_h$, since all other partial derivatives of $J(X,U)$ give the value zero.  For simplicity, we suppress loops over the xyz-indices.

\vspace{.3cm}

\LinesNotNumbered
\begin{algorithm}[H]
\label{algo:dJ_dU}

\KwData{FE global node set $X$, FE connectivity and surface
element tables,\\ FE shape functions $ \hat \theta_j$ and gradients $\widehat{\nabla}\hat\theta_j$,\\
FE solution $U=(u_j)$, \\
FE surface quadrature points and weights,\\
Elasticity and lifing material constants.}

 \KwResult{A  $(3 \times n_{\rm sh} \times N_{\rm fc}^F)$ tensor containing $\frac{\partial J}{\partial U_{\rm loc}}$.}
 
 \vspace{.3cm}
initialization $\frac{\partial J}{\partial U_{\rm loc}}$ $\gets$ 0\;

 \For{  all faces $F\in{\cal N}_h$  }
   {
   \nl initialize a $3\times n_{\rm sh}\times l_q^ F$-tensor $\frac{\partial J}{\partial U}^F$ $\gets$ $0$\;

     \For{  all local degrees of freedom, $k=1,\ldots,n_{\rm sh}$  }
        {

            \For{all surface quadrature points $\hat\xi_l^F$, $l=1,\ldots,l_{q}^F$  }
                {
                  \nl\label{algoline:alg1Start}  Calculate the derivative of \eqref{eqa:Jacobi} and \eqref{eqa:Divergence} with respect to $u_{k}$ at the\\ quadrature point $ \hat\xi_l^F$\;
                    \nl Calculate the derivative of the stress tensor $\sigma(T_{K(F)}(\hat\xi_l^F))$\\ with respect to $u_k$\;
                   \nl Calculate the derivative of $\sigma_v$ with respect to $u_{k}$\;
                    Use \eqref{eqa:Ndet} as well as  \eqref{eqa:RO}, \eqref{eqa:CMB} and \eqref{eqa:ShakeDown} to calculate the $u_{k}$-derivative of ${\rm Ni}_{\rm det}(\sigma(T_{K(F)}(\hat\xi_l^F)))$ \;
                    \nl Use this to calculate the derivative of $\left(\frac{1}{Ni_{\rm det}(\sigma(T_{K(F)}(\hat\xi_l^F)))}\right)^{\bar{m}}$\;
               \nl\label{algoline:alg1End} Multiply this with the surface quadrature weight $\hat\omega^F_l$ and with\\ the Gram determinant $\sqrt{g_F(\hat \xi_l^F)}$\;
                \nl Store the result in $\frac{\partial J}{\partial U}^F[\cdot,k,l]$\;

              }

        }
   \nl    Sum up the $(3 \times n_{\rm sh} \times l_{q}^F)$ tensor $\frac{\partial J}{\partial U}^F $ obtained over the quadratures;\\ 
     \nl  Store the result in $\frac{\partial J}{\partial U_{\rm loc}}[\cdot ,\cdot, F]$ \;
    }

\vspace{.3cm}

 \caption{Compute $\frac{\partial J}{\partial U_{\rm loc}}$ }
\end{algorithm}

\vspace{.3cm}

The previous algorithm provides a surface face-wise separate computation of the partial derivative of the cost functional
$J$ with respect to the local vector of displacements and thus gives a $(3 \times n_{\rm sh} \times N_{\rm fc}^F)$ tensor.
To evaluate the global partial derivative, we have to assemble the tensor
$\frac{\partial J}{\partial U_{\rm loc}}$ based on the surface element and connectivity lists. The following algorithm describes the assembling process.

\begin{algorithm}[H]
\label{algo:dJ_dU_assemble}

\KwData{A $3 \times n_{\rm sh} \times N^F_{\rm fc}$ tensor, FE connectivity table \\
            and surface element table}

\KwResult{A  $3  \times N $ matrix containing $\frac{\partial J}{\partial U}$}

\vspace{.3cm}
initialization $\frac{\partial J}{\partial U}$ $\gets$ 0\;
 \For{  all $F\in\mathcal{N}_h$   }
   {
 
      \For{all local degrees of freedom $k=1,\ldots,n_{\rm sh}$  }
         {
       \nl  $\frac{\partial J}{\partial U}[\cdot,\widehat{j}(K(F),k)]$ $\gets$ $\frac{\partial J}{\partial U}[\cdot,\widehat{j}(K(F),k)]$ + $\frac{\partial J}{\partial U_{\rm loc}}[\cdot,k,F]$;
         }
   }

\vspace{.3cm}

 \caption{Assembling $\frac{\partial J}{\partial U_{\rm loc}}$ to $\frac{\partial J}{\partial U}$ }
\end{algorithm}


The algorithm to calculate $\frac{dJ}{dX}$ is similar to the Algorithm \ref{algo:dJ_dU}, however we have to take the potential dependency of   $\sqrt{g_F(\hat \xi_l^F)}$ and $(\widehat{\nabla} T_{K(F)}(\hat\xi_l))^{-1}$ of the $j$-th node $X_j$ into account. The following algorithm highlights the changes that are needed, detailed calculations are given in Appendix \ref{app:dJ_dX}. Note that certain steps, as e.g. step \ref{algoline:algo3Jacobi} in the following algorithm, can be done beforehand, as the result will only depend on element type and quadrature points. We leave such steps in the loop for better readability.

\vspace{.3cm}

\begin{algorithm}[H]\label{algo:dJ_dX}
\KwData{Same as in Algorithm \ref{algo:dJ_dU}}
 \KwResult{A  $(3 \times n_{\rm sh} \times N_{\rm fc}^F)$ tensor containing $\frac{\partial J}{\partial X_{\rm loc}}$}
 
  \vspace{.3cm}
\nl initialization $\frac{\partial J}{\partial X_{\rm loc}}$ $\gets$ 0\;

  \For{  all faces $F\in{\cal N}_h$  }
   {
    \nl initialize a $3\times n_{\rm sh}\times l_q^F$-tensor $\frac{\partial J}{\partial X}^F$ $\gets$ $0$\;
     \For{  all local degrees of freedom, $k=1,\ldots,n_{\rm sh}$  }
        {
            \For{all surface quadrature points $\hat\xi_l^F$, $l=1,\ldots,l_{q}$  }
                {
                   \nl Calculate the derivative of Jacobian matrix $\widehat{\nabla}T_{K(F)}(\hat\xi_l^F)$\\ with respect to $X_{k}^{K(F)}$\;
                    \nl\label{algoline:algo3Jacobi} Calculate the derivative of  \eqref{eqa:Jacobi} and \eqref{eqa:Divergence}
                    with respect to $X_{k}^{K(F)}$ using \eqref{eqa:Transformation}\;
                    \nl Calculate the $X_{k}$ derivative of $\sqrt{\det g_F(\hat\xi_l^F)}$ using 
                    \eqref{eqa:Gram} and \eqref{eqa:Transformation}\;
                    \nl Multiply the result with $\left(\frac{1}{{\rm Ni}_{\rm det}}\right)^{\bar{m}}$ and $\widehat{\omega}_{lF}$\;
                    \nl Store the result in $\frac{\partial J}{\partial X}^F[\cdot,k,l]$\;
                    \nl Follow the steps \ref{algoline:alg1Start}--\ref{algoline:alg1End} of Algorithm  \ref{algo:dJ_dU} analogously with $u_{k}$ replaced by $X_{k}$\;
                    \nl Augment $\frac{\partial J}{\partial X}^F[\cdot,k,l]$ by the result\;
                }
        }

 \nl    Sum up the $(3 \times n_{\rm sh} \times l_{q}^F)$ tensor $\frac{\partial J}{\partial X}^F $ obtained over the quadratures;\\ 
     \nl  Store the result in $\frac{\partial J}{\partial X_{\rm loc}}[\cdot ,\cdot, F]$ \;

     }

\vspace{.3cm}

\caption{Compute $\frac{\partial J}{\partial X_{\rm loc}}$}
\end{algorithm}
\vspace{.3cm}
As described in the Algorithm \ref{algo:dJ_dU_assemble}, we have now to assemble $\frac{\partial J}{\partial X_{\rm loc}}$ to obtain
$\frac{\partial J}{\partial X}$.




A numerical calculation of $\frac{\partial U(X)}{\partial X}$ e.g. by finite differences is possible in principle, but in many instances is not practically feasible as the dimension of $X$ matches the degrees of freedom, which can be several hundreds of thousands in industrial application. As usually, we avoid the problem of solving a finite element simulation for each of these degrees of freedom by using the Lagrangian formalism in the next subsection.

\subsection{Shape Sensitivity via the  Adjoint Equation}
 We define the Lagrangian \cite{Troe} of the discretized problem as
\begin{equation}
\label{eqa:Lagrange}
\mathscr{L}(X,U,\Lambda)=J(X,U)-\Lambda^T\left(B(X)U-F(X)\right).
\end{equation}
Here, the adjoint state  $\Lambda=(\lambda_j)_{j\in\{1,\ldots,N\}}$ is a $N\times 3$ tensor and the expression $\Lambda^T(\cdots)$ is understood in the sense that the $(j,r)$ indices in \eqref{eqa:WeakPDEDiscreteGlobal} are contracted with the corresponding indices of $\Lambda$. If Dirichlet boundary conditions hold for $U=(u_j)$, we have the same boundary conditions $\lambda_j=0$ for $X_j\in\partial\Omega_D$ for the adjoint state.  As usually, $\frac{\partial\mathscr{L}(X,U,\Lambda)}{\partial \Lambda}=0$ gives the state equation \eqref{eqa:WeakPDEDiscreteGlobal} and
\begin{equation}
\label{eqa:AdjointState}
0=\frac{\partial\mathscr{L}(X,U,\Lambda)}{\partial U}=\frac{\partial J(X,U)}{\partial U}-\Lambda^T B(X) ~~\Leftrightarrow B(X)\Lambda=\frac{\partial J(X,U)}{\partial U}
\end{equation}
is the adjoint equation. Here we used the symmetry of the stiffness matrix in the $(j,r)$ and $(k,l)$ indices, see \eqref{eqa:WeakPDEDiscreteGlobal}. Assuming that the state equation and the adjoint equation hold, the partial derivative of $\mathscr{L}(X,U,\Lambda)$ with respect to $X$ equals the total derivative of $J(X,U(X))$, hence
\begin{equation}
\label{eqa:ShapeGradient}
\frac{dJ(X,U(X))}{dX}=\frac{\partial \mathscr{L}(X,U,\Lambda)}{\partial X}=\frac{\partial J(X,U)}{\partial X}-\Lambda^ T\left(\frac{\partial B(X)}{\partial X}U-\frac{\partial F(X)}{\partial X}\right).
\end{equation}
The above equation is understood in the sense that all tensor indices from $\frac{\partial B(X)}{\partial X}$ and $\frac{\partial F(X)}{\partial X}$ are contracted with $\Lambda$ and $U$, except those that originate from the partial differentiation with respect to $X$. Note that the tensor $\frac{\partial B(X)}{\partial X}$ is of dimension $(N\times 3)\times (N\times 3)\times (N\times 3)$ and  $\frac{\partial F(X)}{\partial X}$ is still of dimension $(N\times 3)\times (N\times 3)$. A one-to-one storage of these quantities in the main memory for real world engineering FE models will be impossible on almost all computer architectures and remain very memory consuming even if sparse data structures are used. We therefore calculate the quantities $\Lambda^T\frac{\partial B}{\partial X}U$ and $\Lambda^T\frac{\partial F}{\partial X}$ directly.
The following algorithm gives the numerical calculation of the shape sensitivity:

\vspace{.3cm}

\begin{algorithm}[H]
\label{algo:ShapeSensitivity}
\KwData{Same as in Algorithm \ref{algo:dJ_dU} plus FE volume quadrature points and weights.}
 \KwResult{A $N\times 3$ tensor containing $\frac{d J}{d X}$}
\nl Use Algorithm \ref{algo:dJ_dU} to obtain $\frac{\partial J}{\partial U}$\;
\nl Solve the adjoint equation \eqref{eqa:AdjointState} numerically using a standard FE solver\;
 \nl Use Algorithm \ref{algo:dB_dX} below and $\Lambda$ to obtain $\Lambda^T\frac{\partial B}{\partial X}U$\;
 \nl Use Algorithm \ref{algo:dFvol_dX} and Algorithm \ref{algo:dFsur_dX} below and $\Lambda$ to obtain $\Lambda^T\frac{\partial F}{\partial X}$\;
\nl Use Algorithm \ref{algo:dJ_dX} to obtain $\frac{\partial J}{\partial X}$\;
\nl Add the tensors from step 4--5 with the proper signs to obtain $\frac{dJ}{dX}$ according to \eqref{eqa:ShapeGradient}.

\vspace{.3cm}

\caption{Compute the Shape Sensitivity $\frac{dJ}{dX}$}
\end{algorithm}

\vspace{.3cm}

The following three algorithms are thus needed for the numerical calculation of the total shape gradient $\frac{dJ}{dX}$. The detailed calculations can be found in Appendix \ref{app:dB_dX} and Appendix \ref{app:dF_dX}. Note that in many cases, it will not be possible to store the $3\times n_{\rm sh}\times 3\times n_{\rm sh}\times  3\times n_{\rm sh}\times N_{\rm el}$ tensor $\frac{\partial B}{\partial X_{\rm loc}}$ due to constraints of the main memory, even if a sparse data format is used. We circumvent this obstacle by contracting with the original solution $U$ and the adjoint state $\Lambda$ during the local calculation. Note that again, for better readability, we did not try to avoid obvious redundant calculations in the loops.

\vspace{.3cm}

\begin{algorithm}[H]

\KwData{Same as in Algorithm \ref{algo:ShapeSensitivity}, but with volume quadrature points and weights.\\
 The adjoint state $\Lambda$.}
 \KwResult{A $ 3\times n_{\rm sh}\times N_{\rm el}$ tensor containing $\Lambda^T\frac{\partial B}{\partial X_{\rm loc}}U$}
 
 \vspace{0.3cm}
 
 \nl initialization $\Lambda^T\frac{\partial B}{\partial X_{\rm loc}}U$ $\gets$ $0$\;

  \For{  all elements $K\in \mathcal{T}_h$  }
        {
        \nl initialize a $3\times n_{\rm sh}\times 3\times n_{\rm sh}\times  3\times n_{\rm sh}\times l_q$ tensor $\frac{\partial B}{\partial X_{\rm loc}}^K$ $\gets$ $0$\;
        \nl initialize a $3\times n_{\rm sh}$ tensor $\Lambda_{\rm loc}^K$ $\gets$ $0$\;
        \nl initialize a $3\times n_{\rm sh}$ tensor $U_{\rm loc}^K$ $\gets$ $0$\;
            \For{  all local degrees of freedom, $j=1,\ldots,n_{\rm sh}$  }
                    {
                    \For{  all local degrees of freedom, $k=1,\ldots,n_{\rm sh}$  }
                    {
                    \For{  all local degrees of freedom, $q=1,\ldots,n_{\rm sh}$  }
                    {
                        \For{all volume quadrature points $\hat\xi_l$, $l=1,\ldots,l_{q}$  }
                            {
                                \nl Compute $\frac{\partial}{\partial X_{j}}\bigl(\omega_{lK}\bigr)$ 
                                using \eqref{eqa:domega_dX} ff.\;
                                \nl Compute $\frac{\partial}{\partial X_{j}}\bigl(\nabla_{[\cdot]} \theta_q(\xi_l)\bigr)$ and $\frac{\partial}{\partial X_{j}}\bigl(\nabla_{[\cdot]} \theta_k(\xi_l)\bigr)$ using \eqref{eqa:dNtheta_dX}\;
                                \nl Compute $\frac{\partial B}{\partial X_{\rm loc}}^K[\cdot,j,\cdot,k,\cdot,q,l]$ using \eqref{eqa:dB_1_dX} and \eqref{eqa:dB_2_dX}\;
                                \nl Assign $\Lambda_{\rm loc}^K[\cdot,k]$ $\gets$ $\Lambda[\cdot,
                                \widehat{j}(K,k)]$\;
	                            \nl Assign $U_{\rm loc}^K[\cdot,q]$ $\gets$ $U[\cdot,
                                \widehat{j}(K,q)]$\;    
                            }
                    }
        }
        \nl Multiply $\frac{\partial B}{\partial X_{\rm loc}}^K[\cdot,j,\cdot,k,\cdot,q,l]$ with $\Lambda_{\rm loc}^K[\cdot,k]$ and $U_{\rm loc}^K[\cdot,q]$ and sum \\ over $q,k$ 
        (along with related $xyz$ indices) and quadrature index $l$\;
        \nl Store the result in  $\Lambda^T\frac{\partial B}{\partial X_{\rm loc}}U[\cdot,j,K]$\;
        }
        }

\caption{\label{algo:dB_dX}Compute $\Lambda^T\frac{\partial B}{\partial X_{\rm loc}}U$ }
\end{algorithm}

\vspace{.3cm}

The assembly of $\Lambda^T\frac{\partial B}{\partial X_{\rm loc}}U$  to   $\Lambda^T\frac{\partial B}{\partial X}U$ is similar to Algorithm \ref{algo:dJ_dU_assemble} and is omitted here. Apparently, Algorithm \ref{algo:dB_dX} will be very run time intensive. However, it can easily be speeded up by embarrassing parallelism over the elements $K\in\mathcal{T}_h$.

The following algorithm gives the derivative of the volume loads with respect to $X$ for the case where the volume force density does not change with $X$. We will discuss the necessary adjustments in the case where volume loads, as e.g. centrifugal loads, depend on $X$ in the context of concrete models in Section \ref{sec:NUM}.


The calculation of the partial derivatives of the volume force term leads to similar storage problems as in Algorithm \ref{algo:dB_dX}. Again we solve this by an element wise calculation and contraction with the adjoint state $\Lambda$. Note that the storage problem is less severe by a $3\times N$ factor for the global expressions. Hence an alternative strategy would be a direct storage of the global data in a sparse matrix format. The force vector might (centrifugal load) or might not (gravity) be $X$-dependent, see the following section.

\vspace{.3cm}

\begin{algorithm}[H]
\label{algo:dFvol_dX}
\KwData{Same as in Algorithm \ref{algo:dJ_dU}, but with volume quadrature points and weights\\
Volume force vector $f$ (eventually depending on $X$)\\
The adjoint state $\Lambda$.}
 \KwResult{A $3\times n_{\rm sh}\times N_{\rm el}$ tensor $\Lambda^T\frac{\partial F^{\rm vol}}{\partial X_{\rm loc}}$}
 \vspace{0.3cm}
 
 \nl initialization $\Lambda^T\frac{\partial F^{\rm vol}}{\partial X_{\rm loc}}$ $\gets$ $0$\;

  \For{  all elements $K\in \mathcal{T}_h$  }
        {
        \nl initialize a $ 3\times n_{\rm sh}\times  3\times n_{\rm sh}\times l_q$ tensor $\frac{\partial F^{K~\rm vol}}{\partial X_{\rm loc}}$ $\gets$ $0$\;
        \nl initialize a $3\times n_{\rm sh}$ tensor $\Lambda_{\rm loc}^K$ $\gets$ $0$\;
            \For{  all local degrees of freedom, $j=1,\ldots,n_{\rm sh}$  }
                    {
                    \For{  all local degrees of freedom, $k=1,\ldots,n_{\rm sh}$  }
                    {
                   
                        \For{all volume quadrature points $\hat\xi_l$, $l=1,\ldots,l_{q}$  }
                            {
                                \nl Compute $\frac{\partial}{\partial X_{j}}\bigl(\omega_{lK}\bigr)$ 
                                using \eqref{eqa:domega_dX} ff.\;
                                \nl Compute $\frac{\partial f_{[\cdot]}}{\partial X_{j}}(\xi_l)$ (model specific)\;
                                \nl Compute $\frac{\partial F^{K~\rm vol}}{\partial X_{\rm loc}^K}[\cdot,j,\cdot,k,l]$ using \eqref{eqa:dVolF_dX} \;
                                \nl Assign $\Lambda_{\rm loc}^K[\cdot,k]$ $\gets$ $\Lambda[\cdot,
                                \widehat{j}(K,k)]$\;    
                            }
                    }
        
        \nl Multiply $\frac{\partial F^{K~\rm vol}}{\partial X_{\rm loc}}[\cdot,j,\cdot,k,\cdot,q,l]$ with $\Lambda_{\rm loc}^K[\cdot,k]$ and sum over $k$ 
        (along with\\ the related $xyz$ index) and quadrature index $l$\;
        \nl Store the result in  $\Lambda^T\frac{\partial F^{\rm vol}}{\partial X_{\rm loc}}[\cdot,j,K]$\;
        }
        }
        
        \caption{Compute $\Lambda^T\frac{\partial F^{\rm vol}}{\partial X_{\rm loc}}$}
\end{algorithm}

\vspace{.3cm}

Also in the case of the $X$-derivative of the surface loads one might adjust the force densities, e.g. in order to keep a total force constant when a surface varies. Again, we refer to Section \ref{sec:NUM} for further discussion.

\vspace{.3cm}

\begin{algorithm}[H]
\label{algo:dFsur_dX}
\KwData{Same as in Algorithm \ref{algo:dJ_dU}\\
Surface force vector $g$ (eventually depending on $X$)\\
The adjoint state $\Lambda$.}
 \KwResult{A $3\times n_{\rm sh}\times N_{\rm fc}^F$ tensor $\Lambda^T\frac{\partial F^ {\rm surf}}{\partial X_{\rm loc}}$}
 \vspace{0.3cm}
 
 \nl initialization $\Lambda^T\frac{\partial F^{\rm surf}}{\partial X_{\rm loc}}$ $\gets$ $0$\;

  \For{  all faces $F\in \mathcal{N}_h$  }
        {
        \nl initialize a $ 3\times n_{\rm sh}\times  3\times n_{\rm sh}\times l_q^F$ tensor $\frac{\partial F^{K ~\rm surf}}{\partial X_{\rm loc}}$ $\gets$ $0$\;
        \nl initialize a $3\times n_{\rm sh}$ tensor $\Lambda_{\rm loc}^{K}$ $\gets$ $0$\;
            \For{  all local degrees of freedom, $j=1,\ldots,n_{\rm sh}$  }
                    {
                    \For{  all local degrees of freedom, $k=1,\ldots,n_{\rm sh}$  }
                    {
                   
                        \For{all surface quadrature points $\hat\xi_l^F$, $l=1,\ldots,l_{q}^F$  }
                            {
                                \nl Compute $\frac{\partial}{\partial X_{j}^{K(F)}}\bigl(\omega_{lF}\bigr)$ 
                                using \eqref{eqa:domega_F_dX} ff.\;
                                \nl Compute $\frac{\partial g_{[\cdot]}}{\partial X_{j}^{K(F)}}(\xi_l)$ (model specific)\;
                                \nl Compute $\frac{\partial F^{K~\rm surf}}{\partial X_{\rm loc}}[\cdot,j,\cdot,k,l]$ using \eqref{eqa:dSurF_dX} \;
                                \nl Assign $\Lambda_{\rm loc}^K[\cdot,k]$ $\gets$ $\Lambda[\cdot,
                                \widehat{j}(K(F),k)]$\;    
                            }
                    }
        
        \nl Multiply $\frac{\partial F^{K~\rm surf}}{\partial X_{\rm loc}}[\cdot,j,\cdot,k,\cdot,q,l]$ with $\Lambda_{\rm loc}^K[\cdot,k]$ and sum over $k$ 
        (along with\\ the related $xyz$ index) and quadrature index $l$\;
        \nl Store the result in  $\Lambda^T\frac{\partial F^{\rm surf}}{\partial X_{\rm loc}}[\cdot,j,F]$\;
        }
        }
        
        \caption{Compute $\Lambda^T\frac{\partial F^{\rm surf}}{\partial X_{\rm loc}}$}
\end{algorithm}

\vspace{.3cm}

The assembly to the global $3\times N$ tensors $\Lambda^T\frac{F^{\rm vol}}{\partial X}$ and  $\Lambda^T\frac{F^{\rm vol}}{\partial X}$ again is standard.
 We have thus provided all the necessary sub-algorithms to Algorithm \ref{algo:ShapeSensitivity}, which allows us an efficient calculation of the shape sensitivity. 

\section{Numerical Examples and Validation}
\label{sec:NUM}
\subsection{Implementation Details}
\label{sec:NUMDet}
 
The implementation of the adjoint method is based on scripts in \texttt{R 3.2.1}  and \texttt{python} as well as the commercial solver \texttt{ABAQUS 6.12}. 
R scripts for the implementation of Algorithm \ref{algo:dB_dX} have been parallelized on up to 12 cores using the \texttt{parallel} package. Local and some of the global tensor calculations are performed on compiled code using tensor summations and contractions provided by the R package \texttt{tensor 1.5}.

Element types available are tetrahedral elements with 4 and 8 DoF and brick elements with 8 or 20 DoF. Quadratures for the surface integrals can be chosen with up to 36 quadrature points, while exact and reduced quadratures are available for volume integrals

Clamped Dirichlet and cyclic boundary conditions have been implemented along with volume forces and surface forces.

\subsection{A Bended Rod under Tensile Loading}
As a first example we consider a situation where one has an intuitive idea where the shape gradient should point to. We thus consider a geometry of a rod which is 6mm long and is bended up to a height of 3mm. Material properties are the same as in Section \ref{sec:DISEx}. The diameter is 1mm.  The boundary condition is Dirichlet $u=0$ at the left and a uniform tensile stress on the right end with $12$ N/mm, which makes up a total force of $18.85$N over the total surface of $1.57$mm${}^2$. The FE model consists out of 6410 nodes and 1302 brick elements of type C3D20R with $n_{\rm sh}=20$ local degrees of freedom. The volume quadrature contains $l_q=8$ points and the surface quadrature contains $l_q^F=36$ quadrature points. The resulting $J$-value is $J=4.20160\times 10^ {-11}
$[cycles${}^{-\bar{m}}]$ which corresponds to a Weibull scale variable of $\eta=154\, 274$ cycles.

In the calculation of the partial derivative $\frac{\partial F(X)}{\partial X}$ in
 \eqref{eqa:ShapeGradient}, we have to take into account the potential $X$ dependence of the physical force density $f$.  In the model considered here, $A$ stands for the portion of the surface where the constant force density is applied in $e_1$ direction.  We adjust the force density $f$ in the way that $f_1(X) A(X)$ is kept constant to the value of $18.85$ N, where $A(X)$ is the area of the surface to which the force is applied after a potential deformation due to a change of $X$. If this is neglected,  the shape gradient would point to the inward direction at $A$  to decrease the surface  and thereby lower to total load. In order to account for effect, notify that $A(X)=\sum_{F\in \mathcal{N}_h, F\subseteq A}\sum_{l=1}^{l_q^F}\omega_{lF}$ and thus an extra term
\begin{equation}
\label{eqa:ForceAdaption}
\frac{\partial f_1(\xi_l^F)}{\partial X_j}=-\frac{f_1}{A(X)^2}\sum_{F\in \mathcal{N}_h, K(F)\in \widehat{j}^{-1}(j)_1, F\subseteq A}\sum_{l=1}^{l_q^F}\frac{\partial}{\partial X_j} \omega_{lF},
\end{equation}
which has to be inserted into the middle term \eqref{eqa:dSurF_dX} in Appendix \ref{app:dF_dX}. $\frac{\partial}{\partial X_j} \omega_{lF}$ is calculated in \eqref{eqa:domega_F_dX}.  Here $\widehat{j}^{-1}(j)_1$ is the projection of $\widehat{j}^{-1}(j)=\{(K_1,n_1),\ldots,(K_k,n_k)\}$ to the first component $\{K_1,\ldots,K_k\}$.

\begin{figure}[t]

\centerline{ \includegraphics[width=.49\textwidth]{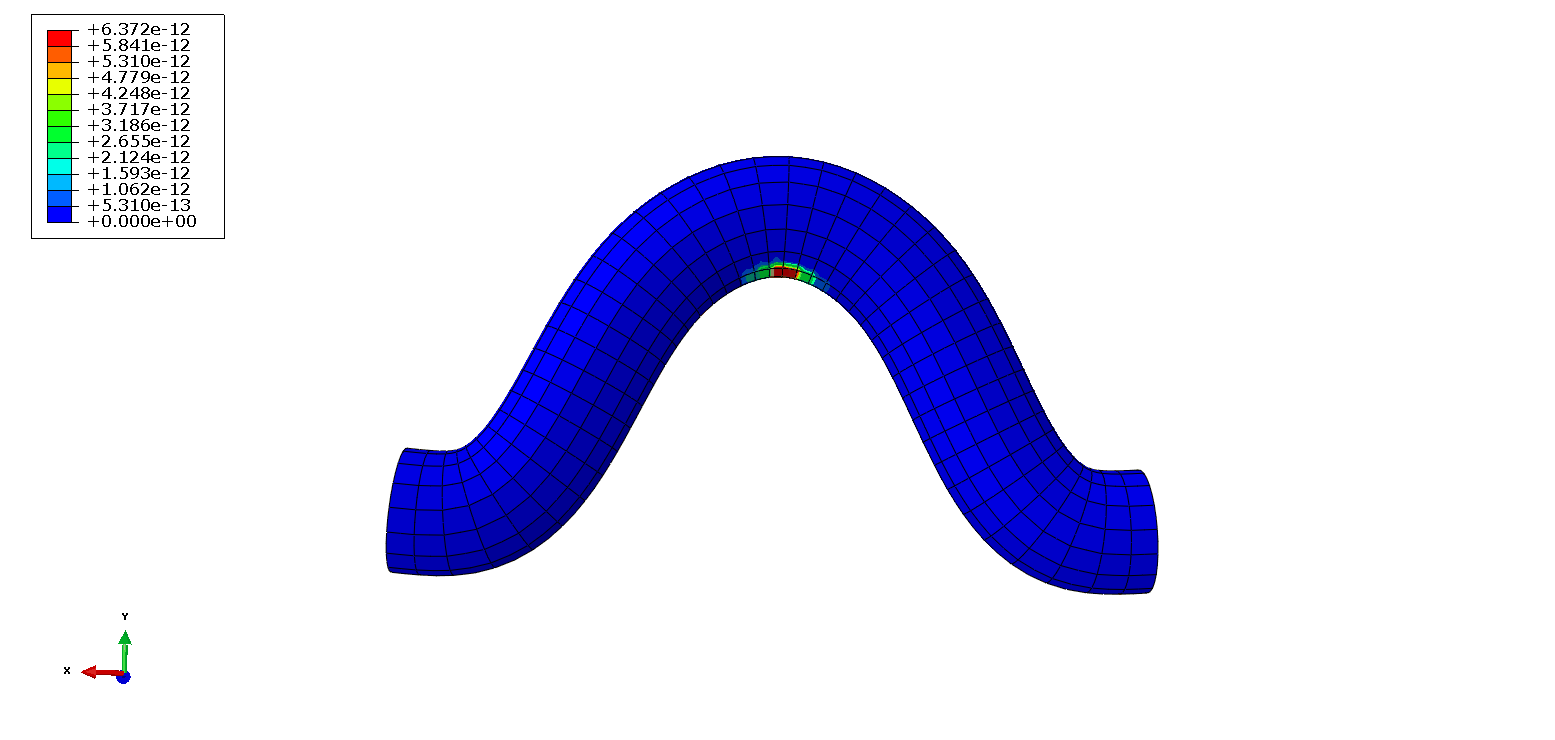}~~~
\includegraphics[width=.49\textwidth]{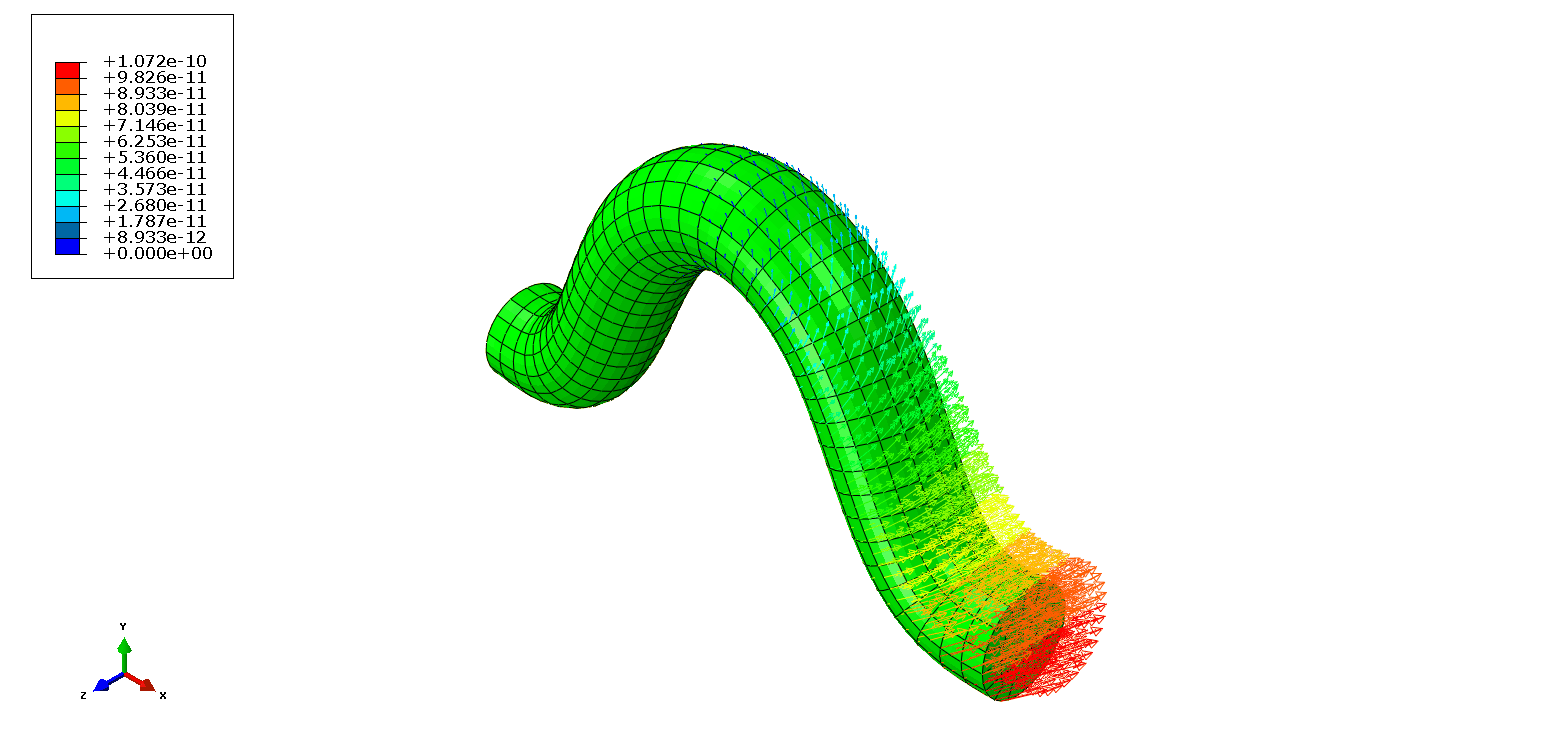}
 }

\caption{\label{fig:IntensityRod} Crack formation intensity and adjoint state: (a) Crack formation intensity for the bended rod. Stress concentration on the lower side of the location with the stongest curvature leads to an augmented probability of stress initiation (left). (b) Adjoint state $\Lambda_j$ is visualized as an arrow at $X_j$. Asymmetry is a consequence of non symmetric boundary conditions.   }
\end{figure}
 
Figure \ref{fig:IntensityRod} (a)  displays the crack formation intensity for the bended rod, whereas Figure \ref{fig:IntensityRod} (b) the adjoint state is displayed. The hot spot for crack initiation is located at the lower side of the portion of the rod with the strongest bending.

\begin{figure}[t]
\centerline{
 \includegraphics[height=.2\textheight]{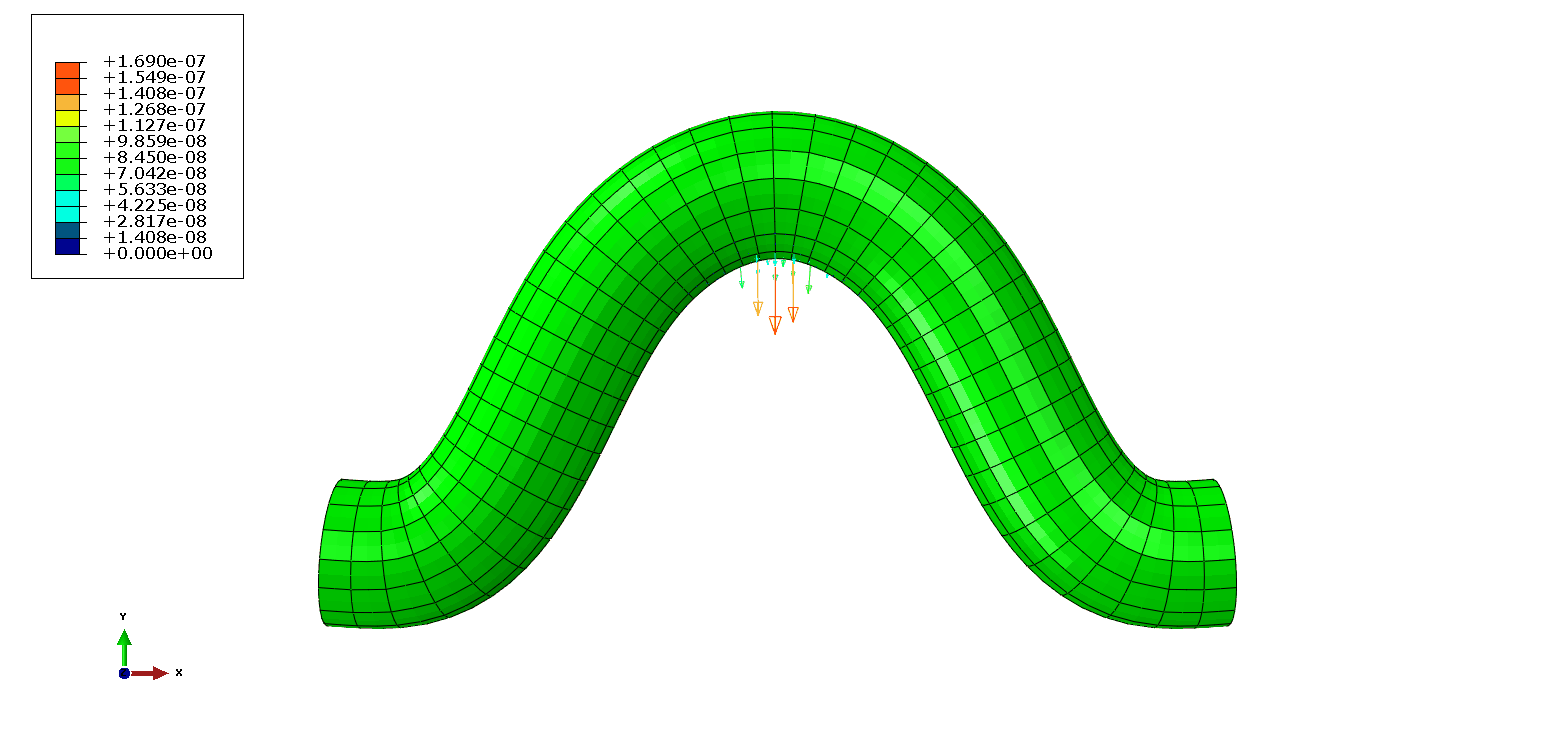}
 \includegraphics[height=.2\textheight]{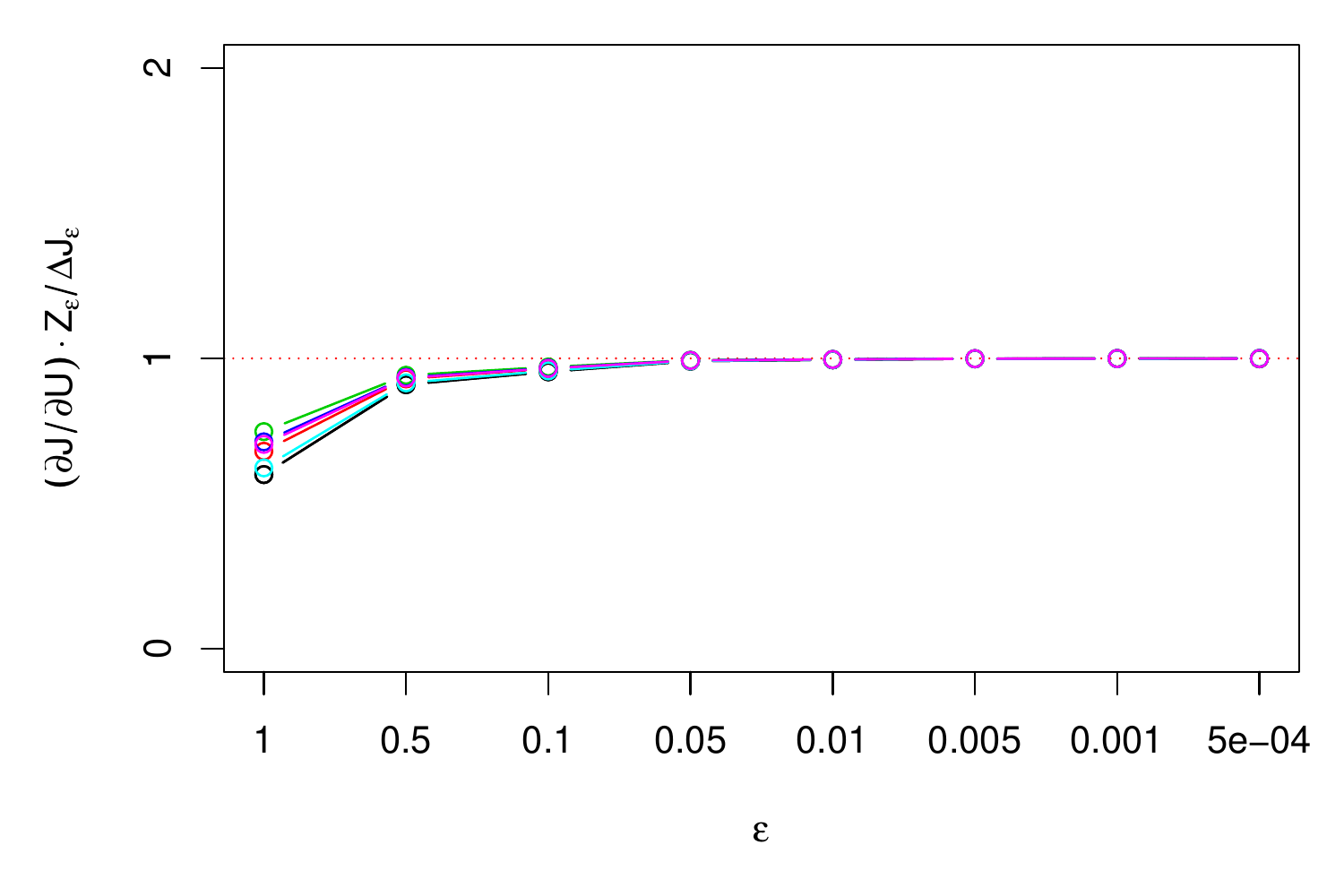}
  }
  
\caption{\label{fig:Jpartial_U} Partial $U$ derivatives of $J$: (a) Partial derivative $\frac{\partial J}{\partial U}_j$ visualized as an arrow at $X_j$ (left). (b)  Comparison of $\frac{\partial J}{\partial U}\cdot V$ with finite the differences $\frac{J(X,U+\varepsilon V)-J(X ,U)}{\varepsilon}$ for $\varepsilon\to 0$. $V$ is taken to be proportional to $U$ (right). }
\end{figure}

 \begin{figure}[t]
\centerline{
 \includegraphics[height=.2\textheight]{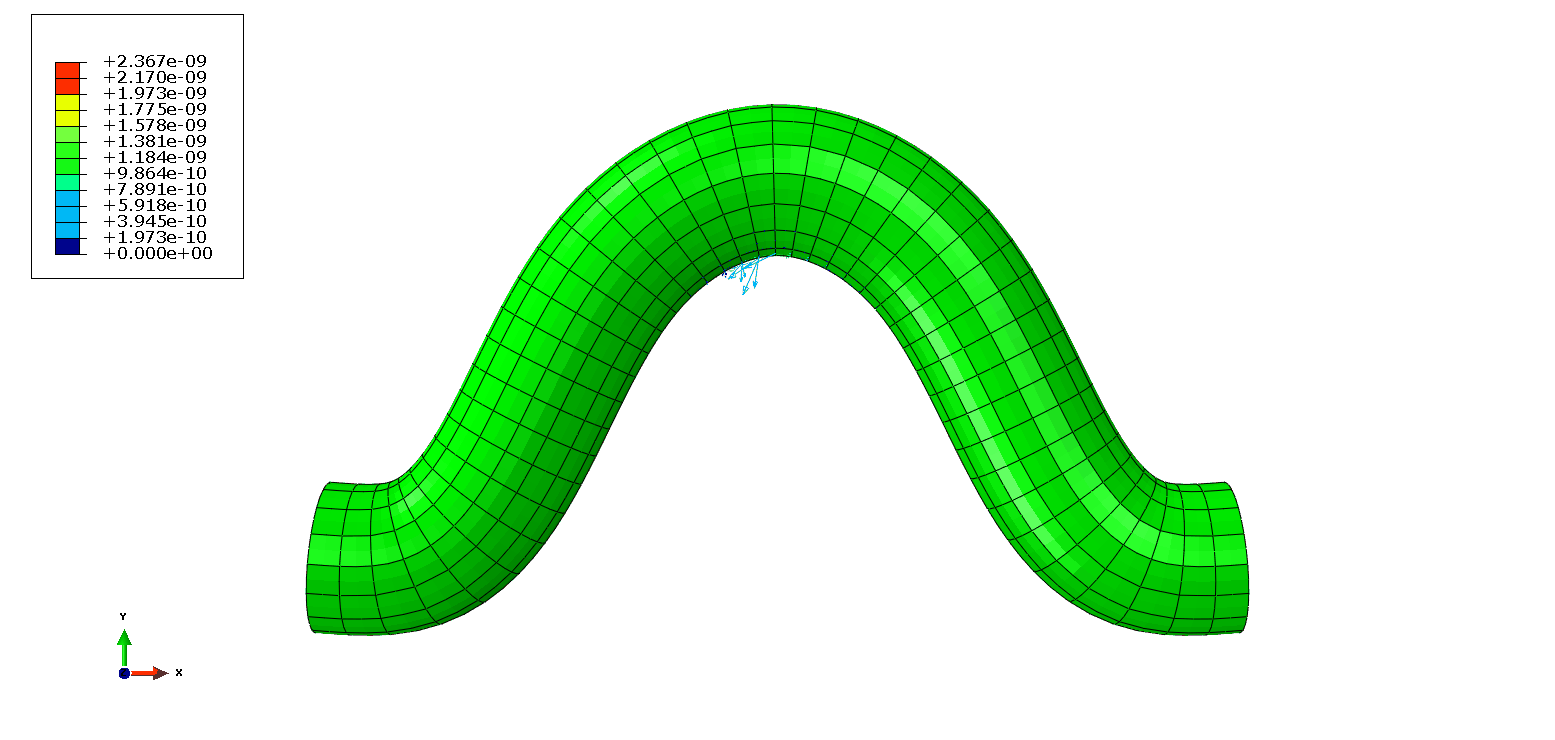}
 \includegraphics[height=.2\textheight]{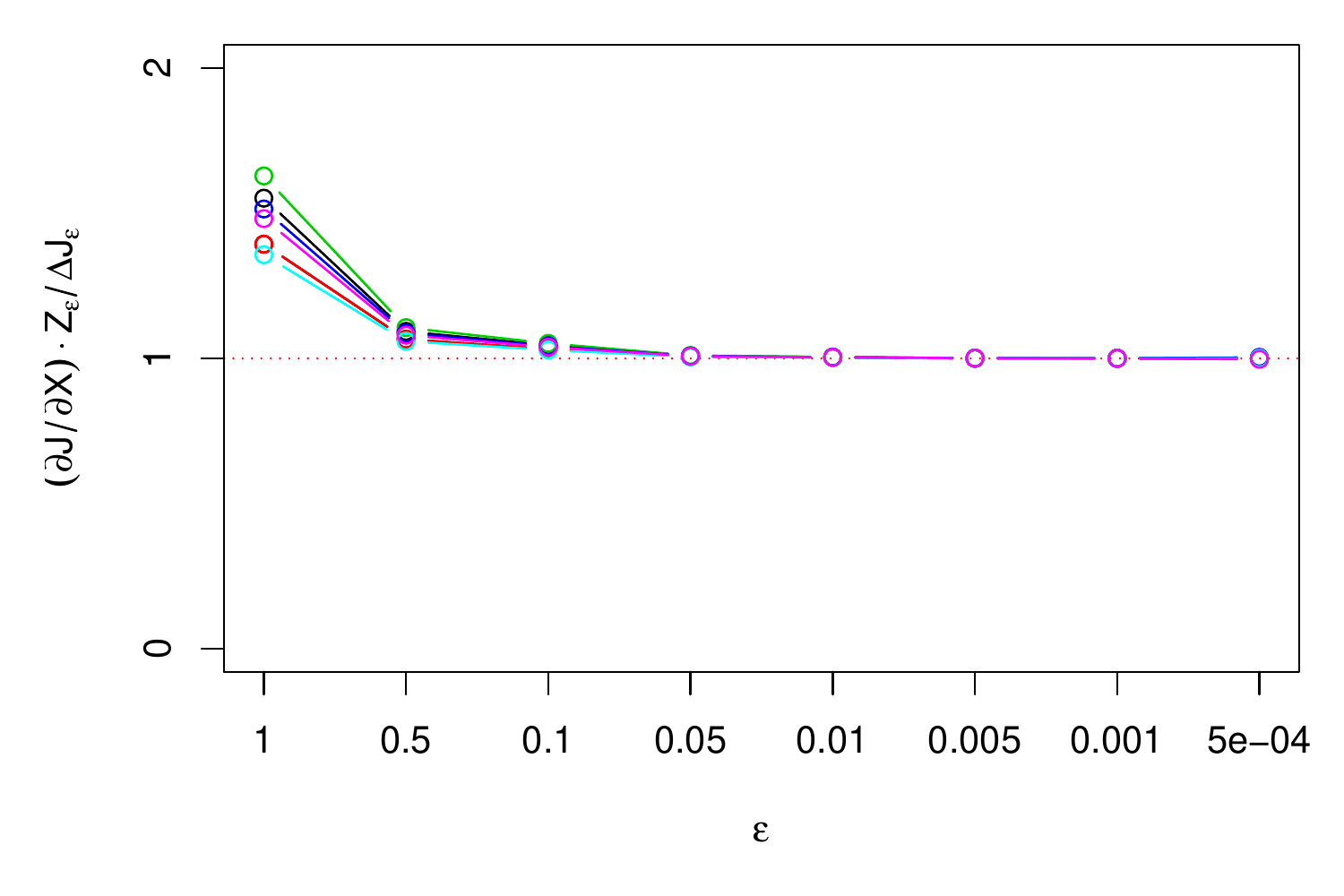}}
  \caption{Partial $X$ derivatives of $J$: (a) The $j$th component of $\frac{\partial J}{\partial X}$ is visualized as an arrow at $X_j$ (left). (b) Comparison of $\frac{\partial J}{\partial X}\cdot V$ with finite the differences $\frac{J(X+\varepsilon V,U)-J(X,U)}{\varepsilon}$ for $\varepsilon\to 0$. $V$ is the normal vector field (right).  }
  \label{fig:Jpartial_X}
\end{figure}

Note that for the adjoint state, Dirichlet conditions are imposed on the left hand side, too, however no such conditions hold on the right. This explains the non symmetric appearance, despite $\frac{\partial J}{\partial U}$ is rather symmetric, see Figure \ref{fig:Jpartial_U} (a).

Figures \ref{fig:Jpartial_U} (a) and \ref{fig:Jpartial_X} (a) display the partial derivative $\frac{\partial J}{\partial U}$ and $\frac{\partial J}{\partial X}$, respectively. The right panels Fig.\ \ref{fig:Jpartial_U} (a)  and Fig.\ \ref{fig:Jpartial_X} (b) compare results of the the finite difference method with the partial derivatives. The numerical validation reveals relative errors in the range of $0.1\%$. 

The components of the derivative $\frac{\partial J}{\partial X_j}$ that are parallel to the surface contour at $X_j$ look artificial. However, such arrows that are pointing away from the location of highest stress can be understood by the interplay of finite element discretization and the probabilistic model. As the function $\left(\frac{1}{{\rm Ni}_{\rm det}}\right)^{\bar{m}}$ is convex in the von Mises stress, the numerics produces a lower integral result for $J$ if an element with high loading is stretched, as this decreases stress values, even if this quenches a neighboring element with potentially lower stress levels. It should however be clear that this is only a by-product of the finite element discretization and not a physical effect. It therefore seems to be reasonable, only to consider normal components of (partial) shape derivatives $\frac{\partial J}{\partial X}$ and $\frac{dJ}{dX}$.

\begin{figure}[ht]

\centerline{ \includegraphics[height=.2\textheight]{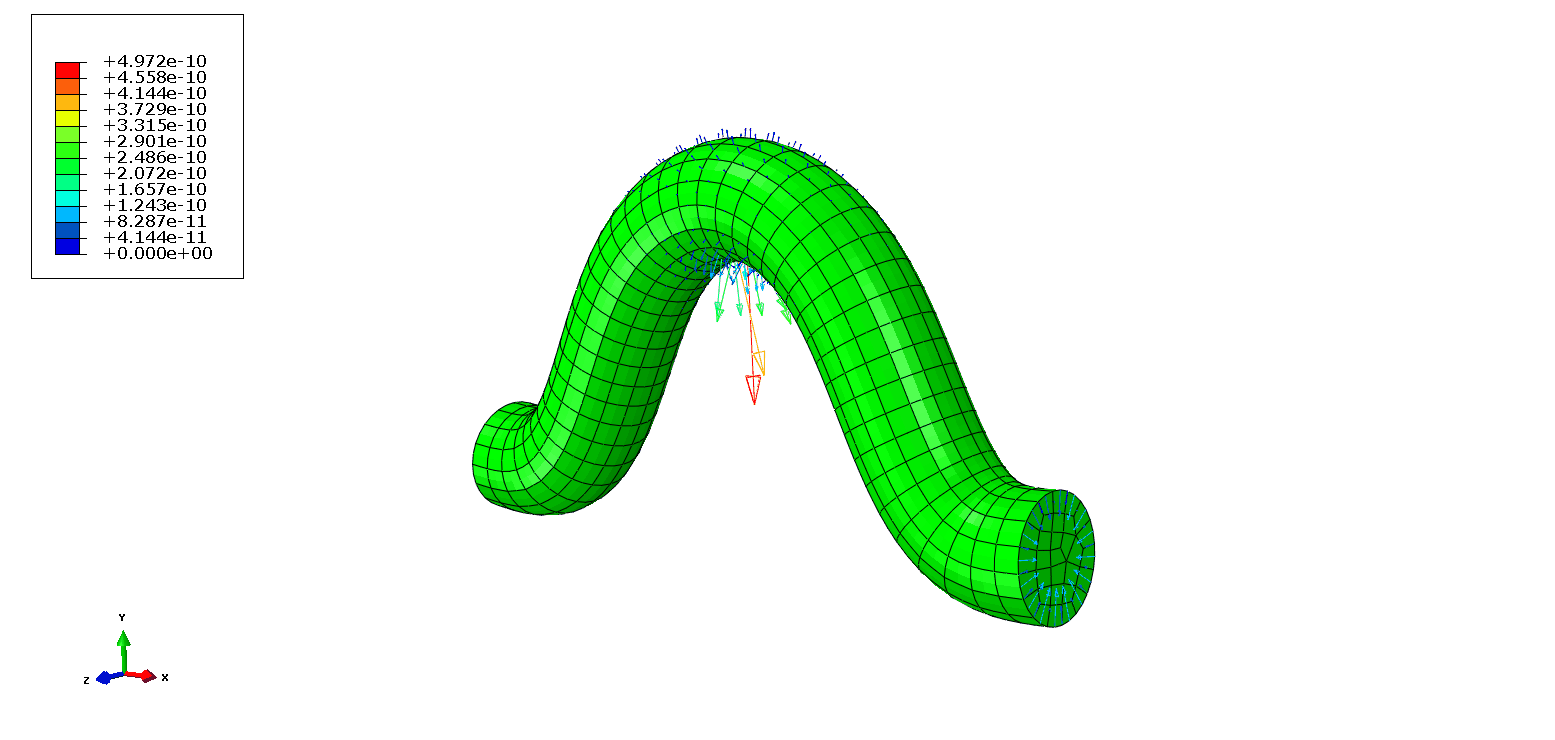}
 \includegraphics[height=.2\textheight]{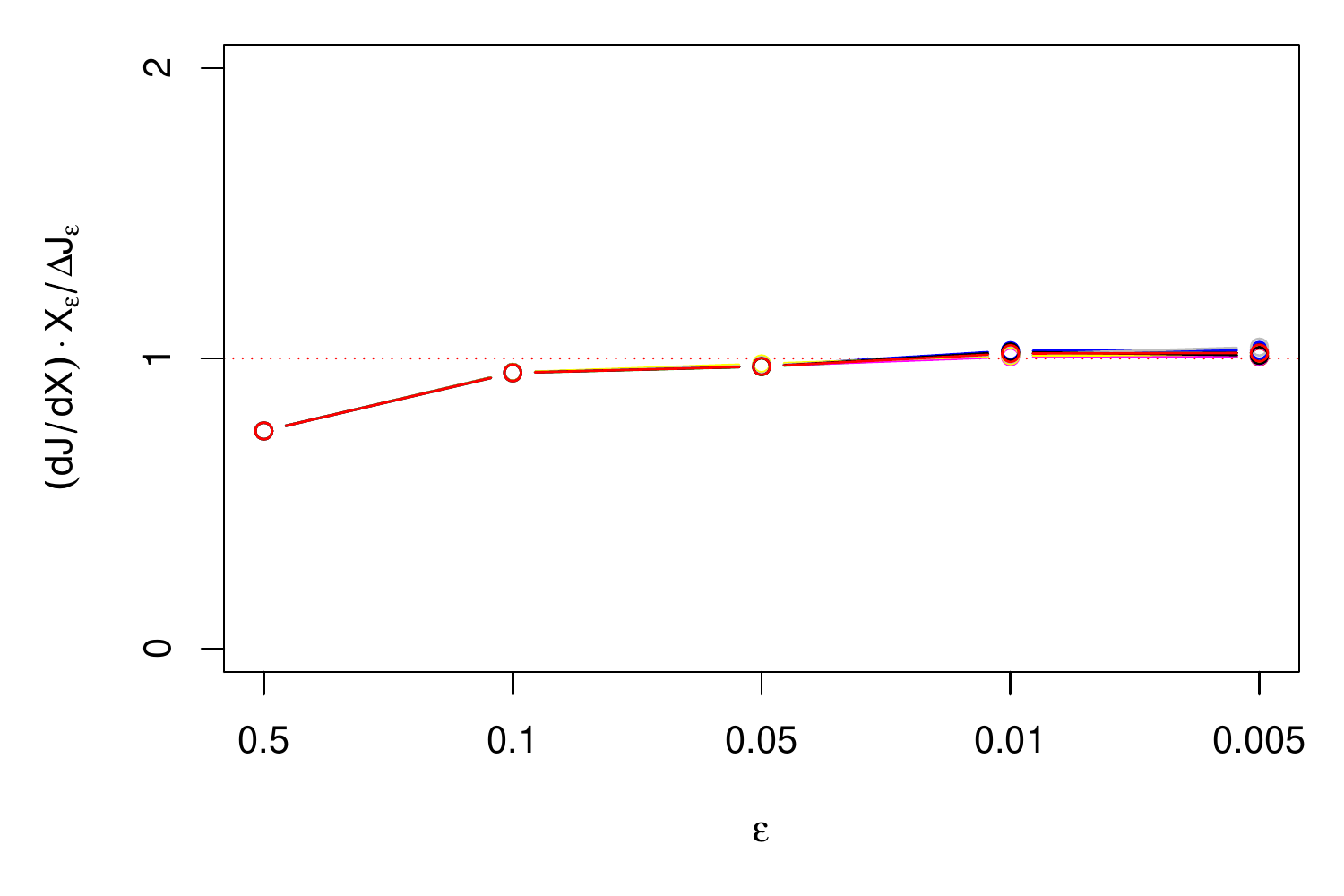}
 }

\caption{\label{fig:sensitivities_rod} Validation of the shape sensitivities: (a) The $j$th component of $\frac{\partial J}{\partial X}$ is visualized as an arrow at $X_j$ (left). (b) Comparison of the shape sensitivity $\frac{d J}{d X}\cdot V$ with finite the differences $\frac{J(X+\varepsilon V,U(X+\varepsilon V))-J(X,U(X))}{\varepsilon}$ for $\varepsilon\to 0$. $V$ is the vector field proportional to $X$ (right). }
\end{figure}

The shape sensitivity is displayed in Figure \ref{fig:sensitivities_rod}. Panel (a)  displays the total shape sensitivity and (b) shows the validation with finite differences in the direction of the node coordinates $X$ which corresponds to a scaling of the structure.


\subsection{Shape Sensitivity of the Radial Compressor}

We next calculate shape sensitivities in the normal direction in the case of the radial compressor model introduced in Section \ref{sec:DISEx}. As gas pressure is neglected in this model, the load vector is exclusively determined by the centrifugal load at the rotation speed of $110 000$ rpm which corresponds to $\omega =2\pi\times 1833.33=11519.17$ Hz and density $\varrho=2650$ kg/m${}^3$. Let $x$ be the rotation axis and $\xi^\perp=(0,\xi_2,\xi_3)$. Then, the centrifugal force density is $f=\rho \omega^ 2\xi^ \perp$ and thus we have to introduce a term that represents the change of the centrifugal load under the change of node coordinates in $\frac{\partial F}{\partial X}$, see Algorithm \ref{algo:dFvol_dX} and \eqref{eqa:dVolF_dX}, as follows 
\begin{equation}
\label{eqa:dfcentrifugal_dX}
\frac{\partial f}{\partial X_{ji}}(\xi_l^K)=\varrho \omega^2 \frac{\partial}{\partial X}_{ji}T_K(\widehat{\xi}_l)^\perp=\begin{cases}0& \mbox{ if }i=1\\
\varrho\omega^2\widehat{\theta}_j(\widehat{\theta}_l)e_i& \mbox{if } i=2,3
\end{cases},
\end{equation}
where $X_{ji}$ is the local node coordinate and we also used \eqref{eqa:Transformation}. 

\begin{figure}[t]

\centerline{ \includegraphics[height=.2\textheight]{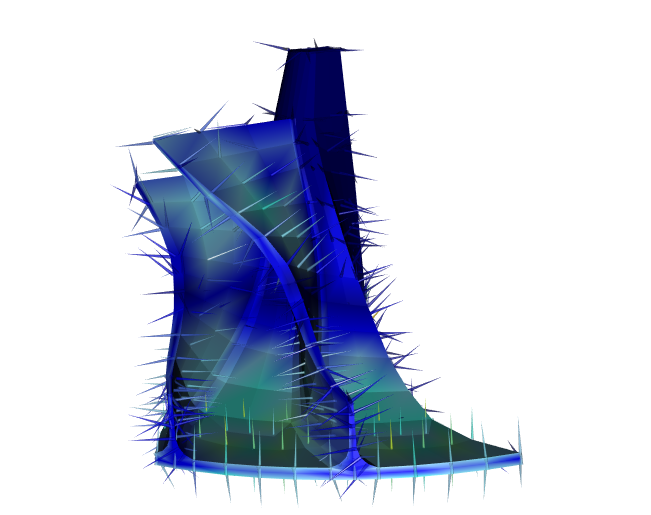}
  \includegraphics[height=.2\textheight]{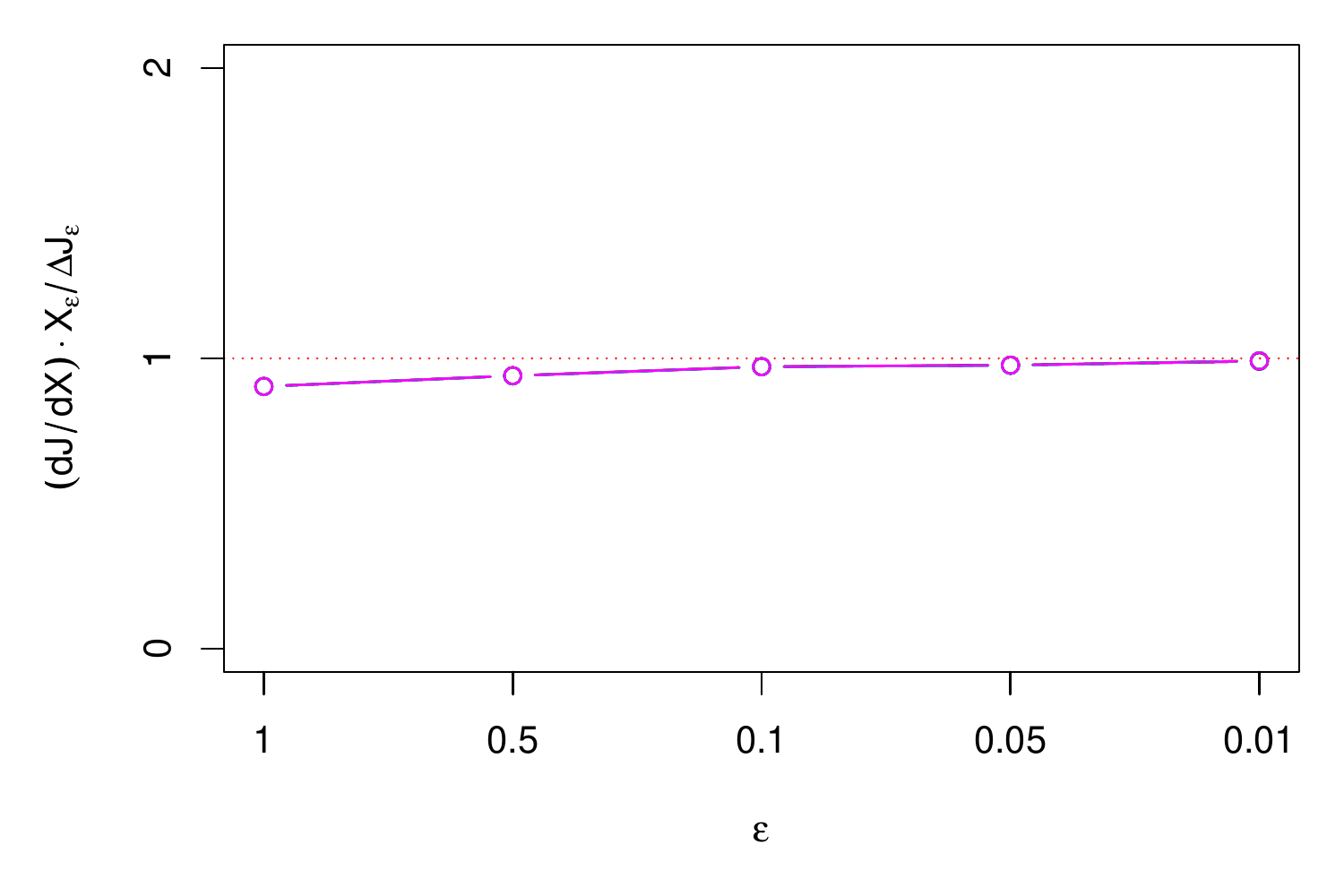}
  }

\caption{\label{fig:SensitivityCompressor} Shape sensitivity for the radial compressor: (a) Shape sensitivity for the radial compressor. Length of $\frac{dJ}{dX}_j$ is displayed as arrow at $X_j$ (left). (b) Validation of the compressor sensitivities vs finite differences. }
\end{figure}

The shape sensitivity for the radial compressor is displayed in figure \ref{fig:SensitivityCompressor}. 
Also in the case of a complex geometry the shape sensitivity calculations conform to more than 99\% with the finite difference method. 
\FloatBarrier
\section{Conclusion and Outlook}
\label{sec:OUT}

In this work we demonstrated that the probability for Low Cycle Fatigue (LCF) crack initiation is an objective functional for mechanical engineering that permits the calculation of shape sensitivities via the discrete adjoint method. We revised the modeling and numerical implementation of probabilistic LCF life models and described the calculation of shape sensitivities for the probability of failure both from the algorithmic and the analytic standpoint. We also provided numerical validation with the relative total error for the shape sensitivities in the range of a few percent, when compared with finite differences. To the best of our knowledge, this is the first numerical implementation of a shape gradient for failure probabilities.  This paves the way for an optimal design for reliability using highly efficient gradient based methods, see \cite{BHGS} for some first steps in 2D. 

We discussed some issues of the discrete adjoint method  related to the enlarging of highly loaded element. This leads to considerable non normal components in the shape derivative. Projection to the normal direction however is  an adequate solution.

It has been demonstrated that the formalism can be implemented with the aid of existing FE solvers using usual desktop architectures with a moderate degree of parallelism. Run times remain moderate even for real world 3D design applications.  Nevertheless, run times can possibly considerably reduced in future implementations, as compared with our R based 'rapid prototype'.

Several directions of  future research open up from this point: From the applied standpoint, the extension of the probabilistic approach to further material classes and damage mechanisms is desirable. Also, having gas turbine design in mind, the treatment of the thermo-mechanical system of partial differential equations is desirable. Furthermore, contact boundary conditions at some point should be integrated into the formalism and the implementation.

It would also be interesting to use the probabilistic objective functional along with the adjoint equations for adaptive mesh refinement as proposed in \cite{BR}.

From a more theoretical standpoint, a better numerical control of the probabilistic objective functional is desirable. Note that this is not straight forward in the present weak and $H^1$-based approach, as no trace on $\partial\Omega$ exists for the stress. Nevertheless, we suggest that with the aid of elliptic regularity theory \cite{AGN,Cia,GS,Schmi} and using derivative recovery schemes, see e.g. \cite{BXZ}, this problem can probably be resolved. 

It also would be desirable to repeat the calculation and implementation of shape derivatives with a first adjoin and then discretize approach and benchmark both strategies for applications in design for optimal reliability. See e.g. \cite{BGS} for some first analytic steps in this direction. This however will not be possible with standard FE solvers and therefore is probably less close to industrial application. A successful benchmark for this strategy could however boost further research and development in that direction with designs with optimal reliability in view.

Last, it is very important to address the challenge of multi objective optimization in a gradient based context. In gas turbine design, shape sensitivities for the aerodynamic performance is already implemented in many publicly available codes. It is therefore highly desirable to put the mechanical design on an equal footing and use the reduction in algorithmic complexity by the adjoint method in multi-physics optimization.

\vspace{.5cm}

\noindent{\bf Acknowledgements:} We thank the Siemens Gas Turbine Engineering Department for Probabilistic Design for valuable discussions and support, in particular thanks go to Dr. Georg Rollmann and Dr. Sebastian Schmitz. We also thank  Dr. T. Seibel (Forschungszentrum J\"ulich) and Prof. Dr. T. Beck (Technical University of Kaiserslautern) and R. Krause (ICS Lugano) for prior joint work on probabilistic LCF, on which this paper is based. This work has been made possible by financial support under the AG Turbo Grant 4.1.13 funded  by Siemens Power and Gas, the German Federal Ministry of Economic Affairs (BMWi) under the grant no. 03ET7041J and by the BMBF project GIVEN under the grant no. 05M2018.

\appendix

\section{Detailed Calculations}
\label{sec:AppCalc}

\subsection{Notation for Finite Elements}
\label{app:FE}
To fix the notation for later purposes, we introduce the standard discretization of the PDE \eqref{eqa:WeakPDE} with Lagrange finite  elements \cite{EG,Bra}.

Let the finite element be given by $\{K,P(K),\Sigma(K)\}$ where $K\subseteq \R^3$ is a compact, connected Lipschitz set, $P$ is a finite vector space  polynomials $p:K\to\R$. The local degrees of freedom are linear functionals on $P(K)$ defined by $\varphi_{K,j}(p)=p(X_j)$ where there are nodes $X_{1}^K,\ldots,X_{n_{\rm sh}}^K\in K$ such that the mapping from $P(K)$ to its local degrees of freedom is bijective. The dual basis of $P(K)$ with respect to the local degrees of freedom are the local shape functions $\theta_{K,k}\in P(K)$ and $\varphi_{K,j}(\theta_{K,k})=\delta_{j,k}$, $j,k\in\{1,\ldots,n_{\rm sh}\}$.

Let $\mathcal{T}_h$ be a finite element mesh on $\Omega$ with $N_{\rm el}$ elements. As usually we assume that there is a reference element $\{\widehat{K},\widehat{P}, \widehat{\Sigma})$ with a finite dimensional linear space of reference polynomials $\widehat{P}$. For each element $K\in\mathcal{T}_h $ in the mesh, we assume that there is a bijective transformation $T_K:\widehat K\to K$ such that $\widehat{P}=P\circ T_K$ , $\widehat\theta_{j}=\theta_{j}\circ T_K$ and $\widehat\varphi_{j}(p\circ T_K)=\varphi_{j}(p)$, $j\in\{1,\ldots,n_{\rm sh}\}$. In many cases, we have
\begin{equation}
\label{eqa:Transformation}
T_K(\hat\xi)=T_K(\hat\xi,X)=\sum_{j=1}^{n_{\rm sh}}\widehat\theta_j(\widehat{\xi})X_{j}^K,~~\widehat\xi\in
\widehat{K}.
\end{equation}

Let $X=\{X_1,...,X_N\}=\bigcup_{K\in \mathcal{T}_h}\{X_{1}^K,...,X_{n_{\rm sh}}^K\}$ be the set of all the Lagrange nodes. For $K\in \mathcal{T}_h$ and $m\in\{1,...,n_{\rm sh}\}$, let $\widehat{j}(K,m)\in\{1,...,N\}$ be the corresponding index of the Lagrange node. The mapping $\widehat{j}(\cdot,\cdot):\mathcal{T}_h\times\{1,\ldots,n_{\rm sh}\}\to\{1,\ldots,N\}$ is also called the connectivity. Let  $\{\theta_{j}:j\in\{1,\ldots,N\}\}$,  be the set of global shape functions $\theta_{j}:\Omega\to\R$. When these functions are restricted to $K\in\mathcal{T}_h$, they belong to $P(K)$ and fulfill
\begin{equation}
\label{eqa:GlShape}
\theta_{j}\restriction_K=\widehat{\theta}_k\circ T_{K}^{-1} \mbox{ provided } \exists k\in\{1,\ldots,n_{\rm sh}\}\mbox{ such that } \widehat{j}(K,k)=j.
\end{equation}
The global finite element space $H_h^1(\Omega,\R)$ is the linear span of $\{\theta_{j}:j\in\{1,\ldots,N\}\}$. Let $H^1_h(\Omega,\R^3)=H_h^1(\Omega,\R)^ {\times 3}$. We have $H_h^1(\Omega,\R^3)\subseteq H^1(\Omega,\R^3)$ and let $H^1_{D,h}(\Omega,\R^3)$ be the subspace of $H^1_h(\Omega,\R^3)$ with functions $u\in H_h^1(\Omega,\R^3)$ vanishing on all boundary nodes $ \overline{\partial\Omega_D}\cap\{X_1,\ldots,X_n\}$. If vanishing of $u$ on the nodes in $\partial\Omega_D$ implies vanishing of $u$ on $\partial\Omega_D$, we have $H_{D,h}^1(\Omega,\R^3)=H_h^1(\Omega,\R^3)\cap H^1_D(\Omega,\R^3)$. A solution $u\in H_{D,h}^1(\Omega,\R^3)$ to the discretized elasticity PDE then fulfills
\begin{equation}
\label{eqa:WeakPDEDiscrete}
   B(u,v)=\int\limits_{\Omega}f\cdot v\, dx+\int\limits_{\partial\Omega_N}g\cdot v \,dA , \forall v\in H_{D,h}^1(\Omega,\R^3).
\end{equation}
Furthermore, such a solution always exists by the coercivity of the bilinear form $B(u,v)$, which also holds on $H^1_{D,h}(\Omega,\R^3)$.

In the derivation of the adjoint equation in Section \ref{sec:AE} we need explicit expressions for both sides of \eqref{eqa:WeakPDEDiscrete} in terms of the global degrees of freedom $U=(u_{j})_{j\in \{1,\ldots,N\}}$, $u_j\in\R^3$, and the node coordinates $X$, where it is understood that $u_{j}=0$ if $X_j\in\partial \Omega_D$.  We rewrite \eqref{eqa:WeakPDEDiscrete} in terms of global variables $U$ via
\begin{equation}
\label{eqa:WeakPDEDiscreteGlobal}
B(X)U=F(X),~~B(X)_{(j,r),(k,s)}=B(e_r\theta_j,e_s\theta_k)\mbox{ and } F_{(j,r)}=\int\limits_{\Omega}f\cdot e_r\theta_j\, dx+\int\limits_{\partial\Omega_N}g\cdot e_r\theta_j \,dA,
\end{equation}
with $e_r$, $r=1,2,3$ the standard Basis on $\R^ 3$. The linear equation \eqref{eqa:WeakPDEDiscreteGlobal} is understood in the sense that the $(k,s)$ indices of the stiffness matrix $B(X)$ are contracted with the related indices of $U$. The solution of this equation is denoted by $U(X)$.

We have emphasised the dependency of the stiffness matrix $B(X)$ and the load vector $F(X)$ from the global node coordinate $N\times 3$-tensor $X$. The explicit calculations, which are completely standard, are given in Appendix \ref{app:DiscretePDE} in order to prepare the ground for the calculation of shape sensitivities.

In some situations, also the surface and volume force densities $f=f(X)$ and $g=g(X)$ are allowed to have a differentiable dependency on $X$. For $f$ this makes sense e.g. in the case of centrifugal loads and we might want to adapt a surface force density $g$ under geometry changes, in order to keep the total force acting on a part of the surface fixed, even if the surface volume changes. For notational simplicity, we will only introduce this when needed for a physically correct description.

\subsection{\label{app:DiscretePDE} Discretization of the PDE}

The integrals on both sides of \eqref{eqa:WeakPDEDiscrete} are implemented  using a numerical quadrature with quadrature points $\widehat{\xi}_l$ and weights $\widehat{\omega}_l$, $l=1,\ldots,l_q$, on the reference element $\widehat{K}$.
We start with the discretization of the left hand side of \eqref{eqa:WeakPDEDiscrete}:
\begin{align}
\label{eqa:discretePDE_A}
\begin{split}
B(u,v)&=\lambda\sum_{K\in \mathcal{T}_h}\int\limits_{K} \nabla\cdot u\nabla\cdot v dx+2\mu\sum_{K\in \mathcal{T}_h}\int\limits_{K}\varepsilon(u):\varepsilon(v) dx\\
&=\lambda\sum_{K\in \mathcal{T}_h}\int\limits_{\widehat{K}} \nabla\cdot u(T_K(\widehat{\xi}))\nabla\cdot v(T_K(\widehat{\xi}))\det(\widehat{\nabla}T_K(\widehat{\xi})) d\widehat{\xi}\\
&~~~~~~~~~~~~~~~~~~~~~+2\mu\sum_{K\in \mathcal{T}_h}\int\limits_{\widehat{K}}\varepsilon\bigl(u(T_K(\widehat{\xi}))\bigr):\varepsilon\bigl(v(T_K(\widehat{\xi}))\bigr)\det(\widehat{\nabla}T_K(\widehat{\xi})) d\widehat{\xi}\\
&=\lambda\sum_{K\in \mathcal{T}_h}\sum_{l=1}^{l_q} \widehat{\omega_l} \det(\widehat{\nabla }T_K(\widehat{\xi_l}))\nabla\cdot u(T_K(\widehat{\xi_l}))\nabla\cdot v(T_K(\widehat{\xi_l}))\\
&~~~~~~~~~~~~~~~~~~~~~+2\mu\sum_{K\in \mathcal{T}_h}\sum_{l=1}^{l_q}\widehat{\omega_l} \det(\widehat{\nabla}T_K(\widehat{\xi_l}))\varepsilon\bigl(u(T_K(\widehat{\xi_l}))\bigr):\varepsilon\bigl(v(T_K(\widehat{\xi_l}))\bigr)
\end{split}
\end{align}
After setting $\omega_{lK}=\widehat{\omega_l} \det \left(\widehat{\nabla} T_K(\widehat{\xi_l})\right)$ and $\xi_l^K=T_K(\widehat{\xi_l})$, this can be written as
\begin{equation}\label{eqa:bilinear}
B(u,v)=\underbrace{\lambda\sum_{K\in \mathcal{T}_h}\sum_{l=1}^{l_q}\omega_{lK}\nabla\cdot u(\xi_l^K)\nabla\cdot v(\xi_l^K)}_{B_1(u,v)}+\underbrace{2\mu\sum_{K\in \mathcal{T}_h}\sum_{l=1}^{l_q}\omega_{lK}\varepsilon\bigl(u(\xi_l^K)\bigr):\varepsilon\bigl(v(\xi_l^K)\bigr)}_{B_2(u,v)}
\end{equation}

If we apply \eqref{eqa:Jacobi} and \eqref{eqa:Divergence} to $u$ and $v$ and insert the result into \eqref{eqa:bilinear}, we obtain the discretization of the left hand side.

For the volume integral on the right hand side of \eqref{eqa:WeakPDEDiscrete} we obtain by a similar argument
\begin{equation}
\label{eqa:DiscreteVolInt}
\int\limits_{\Omega} f\cdot v\, dx=\sum_{K\in\mathcal{T}_h}\sum_{l=1}^{l_q} \omega_{lK} f(\xi_l^K)\cdot v(\xi_l^K).
\end{equation}

And finally we get the following expression for the discretized surface integral
\begin{equation}
\label{eqa:DiscreteSurfInt}
\int\limits_{\partial\Omega}g\cdot v\,dA=\sum_{F\in \mathcal{N}_h}\sum_{l=1}^{l_q^F}\omega_{lF}g(\xi_{lF})\cdot v(\xi_{lF}),\text{where}
~\omega_{lF}=\widehat{\omega}^F_l\sqrt{\det g_F(\widehat{\xi}_l^F)}~\text{and}~\xi_{l}^F=T_{K(F)} (\widehat{\xi}_l^F).
\end{equation}
Here $\widehat{\xi}_l^F$ and $\widehat{\omega}_l^F$ are quadrature points on a reference face $\widehat F$ of the reference element $\widehat{K}$ and $\mathcal{N}_h$ is the collection of all finite element faces that lie in $\partial\Omega$. $g_F(\hat\xi)=\widehat{\nabla}^{\widehat{F}}(T_K\restriction \widehat{F})(\xi)\widehat{\nabla}^{\widehat{F}}(T_K\restriction \widehat{F})^T(\xi)$ is the Gram matrix on $\widehat{F}$, where $T_K:\widehat{F}\to F$.

If the selected quadratures are not exact, the above identities have to be understood in the sense of approximations.

\subsection{\label{app:dJ_dU}Computing $\frac{\partial J}{\partial U}$}
In the following, we use some conventions related to the connectivity mapping $\widehat{j}:\mathcal{T}_h\times\{1,\ldots,n_{\rm sh})\to\{1,\ldots,N\}$. For $K\in \mathcal{T}_h$, we denote by $\widehat{j}_K:\{1,\ldots,n_{\rm sh}\}\to \{1,\ldots,N\}$ the restriction of the connectivity mapping $\widehat{j}$ to the set $\{(K,1),\ldots,(K,n_k)\}$, where we identify $1,\ldots,n_{\rm sh}$ with $(K,1),\ldots,(K,{n_{\rm sh}})\}$.

For a fixed global index $j\in \{1,2,...,N\}$ we have $\widehat{j}^{-1}(j)=\{(K_1,m_1),...,(K_f,m_f)\}$, where $K_1,...,K_f \in \mathcal{T}_h$ and $m_1,...,m_f\in \{1,2,...,n_{\rm sh}\}$. With $\widehat{j}^{-1}(j)_1=\{K_1,\ldots,K_f\}$ we denote the set projection to the first component. 

For $k=1,2,3$ let $u_{jk}$ be the $x,y,z$ coordinate of the global degree of freedom $u_j$. We obtain for the partial derivative of $J$ with respect to this variable
\begin{align}
\label{eqa:dJ_du}
\begin{split}
\frac{\partial J(X,U)}{\partial u_{jk}} &= \sum_{F\in\mathcal{N}_h\atop j\in \widehat{j}_{K(F)}(\{1,\ldots,n_{\rm sh}\})}\sum_{l=1}^{l_q^F}\omega_{lF}\\&\times \frac{\partial}{\partial u_{jk}}\left[{\rm CMB}^{-1}\left(\varepsilon_{a}^{\rm el-pl}\biggl(\sum_{m=1}^{n_{\rm sh}}u_{\widehat{j}(K(F),m)}\otimes((\widehat{\nabla} T_{K(F)}^{T}(\widehat{\xi}_l^F))^{-1}\widehat{\nabla} \widehat{\theta}_m(\widehat{\xi}_l^F))\biggr)\right)\right]^{-\bar{m}}.
\end{split}
\end{align}
We will frequently use the abbreviation
\begin{equation}
\label{eqa:q} 
q=q_{lF}=q_{lF}(U,X)=\nabla u(\xi_l^F)=\sum_{m=1}^{n_{\rm sh}}u_{j_{(K(F),m)}}\otimes (\widehat{\nabla} T_{K(F)}^{T}(\widehat{\xi}_l^F))^{-1}\widehat{\nabla} \widehat{\theta}_m(\widehat{\xi_l})), ~l=1,\ldots,l_q^ F,
\end{equation} 
with $\xi_l^F=T_{K_F}(\hat\xi_l)$. We thus have
\begin{align}
\frac{\partial}{\partial u_{jk}}[{\rm CMB}^{-1}(\varepsilon_{a}^{\rm el-pl}(q_{lF}))]^{-\bar{m}}&=-\bar{m}[{\rm CMB}^{-1}(\varepsilon_{a}^{\rm el-pl}(q_{lF}))]^{-\bar{m}-1}\frac{\partial}{\partial u_{jk}}[{\rm CMB}^{-1}(\varepsilon_{a}^{\rm el-pl}(q_{lF}))].
\end{align}
In the next step, we compute
\begin{align}
\frac{\partial}{\partial u_{jk}}[{\rm CMB}^{-1}(\varepsilon_{a}^{\rm el-pl}(q_{lF}))] &= \frac{\frac{\partial}{\partial u_{jk}}\bigl(\varepsilon_{a}^{\rm el-pl}(q_{lF})\bigr)}{\frac{\partial {\rm CMB}}{\partial \varepsilon_a^{\rm el-pl}}({\rm CMB}^{-1}(\varepsilon_{a}^{\rm el-pl}(q_{lF})))}
\end{align}
We  now calculate 
\begin{align}
  \frac{\partial}{\partial u_{jk}}\bigl(\varepsilon_{a}^{\rm el-pl}(q_{lF})\bigr) = \frac{\partial \varepsilon_{a}^{\rm el-pl}(q_{lF})}{\partial q_{lF}}:\frac{\partial q_{lF}}{\partial u_{jk}}= \tr \biggl(\Bigl(\frac{\partial \varepsilon_{a}^{\rm el-pl}(q_{lF})}{\partial q_{lF}} \Bigr)^T\frac{\partial q_{lF}}{\partial u_{jk}}\biggr),
\end{align}
where for $s,n=1,2,3$,
\begin{align}
\begin{split}
\frac{\partial (q_{lF})_{sn}}{\partial u_{kj}}&=\frac{\partial }{\partial u_{jk}}\left(\sum_{m=1}^{n_{\rm sh}}u_{\widehat{j}(K(F),m)s}\left((\widehat{\nabla}T_{K(F)}^{T}(\widehat{\xi}_l^F))^{-1}\widehat{\nabla} \widehat{\theta}_m(\widehat{\xi}_l^F)\right)_n\right)\\
                     &=\sum_{m=1}^{n_{\rm sh}}\frac{\partial u_{\widehat{j}(K(F),m)s}}{\partial u_{jk}}\left((\widehat{\nabla}T_{K(F)}^{T}(\widehat{\xi}_l^F))^{-1}\widehat{\nabla} \widehat{\theta}_m(\widehat{\xi}_l^F)\right)_n\\
                     &=\begin{cases} \delta_{sk}\left((\widehat{\nabla}T_{K(F)}^{T}(\widehat{\xi}_l^F))^{-1}\widehat{\nabla} \widehat{\theta}_{\widehat{j}^{-1}_{K(F)}(j)}(\widehat{\xi}_l^F)\right)_n& \mbox{ if }j\in \widehat{j}_{K(F)}(\{1,\ldots,n_{\rm sh}\})\\0&\mbox{ otherwise}\end{cases},
\end{split}
\end{align}
with $\delta_{sk}$ Kronecker's symbol.
In order to simplify the calculation, we square \eqref{eqa:ShakeDown} and divide by $E$. The resulting relation between $\varepsilon_{a}^{\rm el-pl}$ and $\frac{\sigma_a^2}{E}$ then is denoted by $\varepsilon_{a}^{\rm el-pl}=\bar{{\rm SD}}\left(\frac{\sigma_a^2}{E}\right)$. The partial derivative of $\varepsilon_{a}^{\rm el-pl}$ with respect to $q_{sm}$ is now:
\begin{align}
\frac{\partial\varepsilon_{a}^{\rm el-pl}}{\partial q_{sm}}=\frac{d {\rm RO}\circ \bar{\rm SD}^{-1}\left(\frac{(\sigma_a)^2}{E}\right)}{d \left(\frac{(\sigma_a)^2}{E}\right)}\cdot\frac{\partial}{\partial q_{sm}}\biggl(\frac{(\sigma_a)^2}{E}\biggr).
\end{align}
The derivative of the function ${\rm RO}\circ\bar{\rm SD}^{-1}(\cdot)$ is easily calculated from \eqref{eqa:RO} and \eqref{eqa:ShakeDown}.
We further calculate
\begin{align}
\label{eqa:dssq}
\frac{\partial}{\partial q_{sm}}\biggl(\frac{(\sigma_a)^2}{E}\biggr)=\frac{3\mu^2}{4E}\sum_{i,j=1}^{3}A_{ij}\frac{\partial A_{ij}}{\partial q_{sm}}
\end{align}
using $A_{ij}=\biggl(-\frac{2}{3}\bigl(\frac{1}{2}(q_{ij}+q_{ji})\bigr)\delta_{ij}+q_{ij}+q_{ji}\biggr)$. We thus obtain
\begin{align}
\frac{\partial A_{ij}}{\partial q_{s,m}}=-\frac{2}{3}\delta_{ij}\delta_{sm}+\delta_{is}\delta_{jm}+\delta_{js}\delta_{im}.
\end{align}

We thus get for the left hand side of \eqref{eqa:dssq}
\begin{align}
\frac{3\mu^2}{4E}\sum_{i,j=1}^{3}\biggl(-\frac{2}{3}\bigl(\frac{1}{2}(q_{ij}+q_{ji})\bigr)\delta_{ij}+q_{ij}+q_{ji}\biggr)\bigl(-\frac{2}{3}\delta_{ij}\delta_{sm}+\delta_{is}\delta_{jm}+\delta_{js}\delta_{im}\bigr).
\end{align}

\subsection{\label{app:dJ_dX}Computing $\frac{\partial J}{\partial X}$}
As in Subsection \ref{app:dJ_dU}, We use the abbreviation $q_{lF}$ for $\nabla u(\xi_l^F)$, see \eqref{eqa:q}. The partial derivative of $J(X,U)$ w.r.t. the global $j$th mesh node $i$-coordinate $X_{ji}$, $i=1,2,3$ and $j=1,\ldots,N$, is
\begin{align}
\label{eqa:dJ_dX}
\begin{split}
\frac{\partial J(X,U)}{\partial  X_{ji}}&=\sum_{F\in\mathcal{N}_h\atop j\in \widehat{j}_{K(F)}(\{1,\ldots,n_{\rm sh}\})}\sum_{l=1}^{l_q^F}\frac{\partial}{\partial  X_{ji}}\Biggl(\omega_{lF}\left({\rm CMB}^{-1}\left(\varepsilon_{a}^{\rm el-pl}(q_{lF})\right)\right)^{-\bar{m}}\Biggr)\\
                                             &=\sum_{F\in\mathcal{N}_h\atop j\in \widehat{j}_{K(F)}(\{1,\ldots,n_{\rm sh}\})}\sum_{l=1}^{l^F_q}\Biggl[\frac{\partial \omega_{lF}}{\partial  X_{ji}}\Bigl({\rm CMB}^{-1}(\varepsilon_{a}^{\rm el-pl}(q_{lF})))\Bigr)^{-\bar{m}}\\
&+ \omega_{lF}\frac{\partial}{\partial  X_{ji}}\Bigl(\Bigl({\rm CMB}^{-1}(\varepsilon_{a}^{\rm el-pl}(q_{lF}))\Bigr)^{-\bar{m}}\Bigr)\Biggr].
\end{split}
\end{align}
We compute first the partial derivative of $\frac{\partial \omega_{lF}}{\partial  X_{ji}}$ as
\begin{align}
\label{eqa:domega_F_dX}
\begin{split}
\frac{\partial \omega_{lF}}{\partial  X_{ji}} &= \frac{\partial}{\partial  X_{ji}}\left(\widehat{\omega_l}\sqrt{\det(g_F(\widehat{\xi}_l^F))}\right)\\
                                              &=\widehat{\omega_l}\frac{1}{2}\bigl(\det(g_F(\widehat{\xi}_l^F))\bigr)^{-1/2}\frac{\partial}{\partial X_{ji}}(\det(g_F(\widehat{\xi}^F_l))).
\end{split}
\end{align}
The derivative of the determinant is
\begin{align}
\begin{split}
\frac{\partial}{\partial X_{ji}}\biggl(\det\Bigl(g_F(\widehat{\xi}_l^F)\Bigr)\biggr)=\det\Bigl(g_F(\widehat{\xi}^F_l)\Bigr)
\tr\biggl(\Bigl(g_F(\widehat{\xi}^F_l)\Bigr)^{-1}\frac{\partial g_F(\widehat{\xi}^F_l)}{\partial X_{ji}}\biggr),
\end{split}
\end{align}
where $g_F(\widehat{\xi}_l^F)=J_F(\widehat{\xi}_l^F)^T J_F(\widehat{\xi}_l^F)$ and
\begin{equation}
J_F(\widehat{\xi}_l^F)=\frac{\partial T_K(\widehat{\xi}_l^F)}{\partial \widehat{X}^F}
=\sum_{r=1}^{n_{\rm sh}} \left(
  \begin{array}{cc}
    X_{r1}^{K(F)}\frac{\partial \widehat{\theta}_r(\widehat{\xi}_l^F)}{\partial \widehat{X}_1^F} & X_{1r}^{K(F)}\frac{\partial \widehat{\theta}_r(\widehat{\xi}_l^F)}{\partial \widehat{X}_2^F} \\
    X_{r2}^{K(F)}\frac{\partial \widehat{\theta}_r(\widehat{\xi}_l^F)}{\partial \widehat{X}_1^F} & X_{r2}^{K(F)}\frac{\partial \theta_r(\widehat{\xi}_l^F)}{\partial \widehat{X}^F_2} \\
    X_{r3}^{K(F)}\frac{\partial \widehat{\theta}_r(\widehat{\xi}_l^F)}{\partial \widehat{X}_1^F} & X_{r3}^{K(F)}\frac{\partial \widehat{\theta}_r(\widehat{\xi}_l^F)}{\partial \widehat{X}^F_2} \\
  \end{array}
\right).
\end{equation}
Here $\widehat{X}^F_{i}$, $i=1,2$, are the coordinates on the two dimensional reference face $\widehat{F}$ corresponding to $F$ in $\widehat{K}$. 

The derivative $\frac{\partial g_F(\widehat{\xi}_l^F)}{\partial X_{ji}}$ is thus given by
\begin{equation}
\frac{\partial g_F(\widehat{\xi}_l^F)}{\partial X_{ji}}=\frac{\partial}{\partial X_{ji}}\bigl(J_F(\widehat{\xi}_l^F)\bigr)^T J_F(\widehat{\xi}_l^F)+\bigl(J_F(\widehat{\xi}_l^F)\bigr)^T\frac{\partial}{\partial X_{ji}}\bigl(J_F(\widehat{\xi}_l^F)\bigr).
\end{equation}
 Furthermore, for $s=1,2,3$ and $k=1,2$ we have
\begin{align}
\begin{split}
\frac{\partial}{\partial X_{ji}}\bigl(J_F(\widehat{\xi}_l^F)_{sk}\bigr)&= \sum_{r=1}^{n_{\rm sh}}\frac{\partial}{\partial X_{ji}}\left(X_{rs}^{K(F)}\frac{\partial \theta_r(\widehat{\xi}_l^F)}{\partial\widehat{X_k}} \right)\\
& =\begin{cases}
  \delta_{is}\frac{\partial \widehat{\theta}_{\widehat{j}^{-1}_{K(F)}(j)}(\widehat{\xi}_l^F)}{\partial\widehat{X_k}}&\mbox{ if } j\in \widehat {j}_{K(F)}(\{1,\ldots,n_{\rm sh}\})\\0&\mbox{ otherwise}
  \end{cases}.
\end{split}
\end{align}
This finishes the computation of the first term on the right hand side of \eqref{eqa:dJ_dX}.

To compute the second term in  \eqref{eqa:dJ_dX}, we take the partial derivative
\begin{align}
\begin{split}
\frac{\partial}{\partial  X_{ji}}\Bigl[{\rm CMB}^{-1}(\varepsilon_{a}^{\rm el-pl}(q_{lF}))\Bigr]^{-\bar{m}}&= \frac{-\bar{m}[{\rm CMB}^{-1}(\varepsilon_{a}^{\rm el-pl}(q_{lF}))]^{-\bar{m}-1}}{\frac{\partial {\rm CMB}}{\partial \varepsilon_{a}^{\rm el-pl}}({\rm CMB}^{-1}(\varepsilon_{a}^{\rm el-pl}(q_{lF})))}
\frac{\partial}{\partial X_{ji}}\bigl(\varepsilon_{a}^{\rm el-pl}(q_{lF})\bigr),
\end{split}
\end{align}
with
\begin{align}
  \frac{\partial}{\partial X_{ji}}\bigl(\varepsilon_{a}^{\rm el-pl}(q_{lF})\bigr) &= \frac{\partial \varepsilon_{a}^{\rm el-pl}(q_{lF})}{\partial q_{lF}}:\frac{\partial q_{lF}}{\partial X_{ji}}.
\end{align}
In the above equation, $:$ stands for the contraction of both $q$ indices. 
Next we have to compute 
\begin{align}
\label{eqa:dq_dX}
\begin{split}
\frac{\partial q_{lF}}{\partial X_{ji}}(X)&=\frac{\partial }{\partial X_{ji}}\biggl(\sum_{m=1}^{n_{\rm sh}}u_{\widehat{j}(K(F),m)}\otimes((\widehat{\nabla}T_{K(F)}^T(\widehat{\xi}_l^F))^{-1}\widehat{\nabla} \widehat{\theta}_m(\widehat{\xi}_l^F))\biggr)\\
                                  &=\sum_{m=1}^{n_{\rm sh}}u_{j(K(F),m)}\otimes\Bigl(\frac{\partial }{\partial X_{ji}}(\widehat{\nabla}T_{K(F)}^T(\widehat{\xi}_l^F))^{-1}\Bigr)\widehat{\nabla} \widehat{\theta}_m(\widehat{\xi}_l^F),
\end{split}
\end{align}
where
\begin{equation}
\label{eqa:dNT_dX}
\frac{\partial}{\partial X_{ji}}\bigl(\widehat{\nabla} T_{K}^{T}(\widehat{\xi}_l^ F)\bigr)^{-1}=-\bigl(\widehat{\nabla} T_{K}^T(\widehat{\xi}_l^ F)\bigr)^{-1}\frac{\partial \bigl(\widehat{\nabla} T_{K}^T(\widehat{\xi}_l^ F)\bigr)}{\partial X_{ji}}\bigl(\widehat{\nabla} T_{K}^T(\widehat{\xi}_l^ F)\bigr)^{-1}.
\end{equation}
The Jacobian matrix has the form
\begin{equation}
\widehat{\nabla} T_{K}(\widehat{\xi}_l^ F)=
\left(
  \begin{array}{ccc}
    \frac{\partial X_{1}^{K}}{\partial \widehat{X_1}} & \frac{\partial X_{1}^{K}}{\partial \widehat{X_2}} & \frac{\partial X_{1}^{K}}{\partial \widehat{X_3}} \\
    \frac{\partial X_{2}^{K}}{\partial \widehat{X_1}} & \frac{\partial X_{2}^{K}}{\partial \widehat{X_2}} & \frac{\partial X_{2}^{K}}{\partial \widehat{X_3}}\\
    \frac{\partial X_{3}^{K}}{\partial \widehat{X_1}} & \frac{\partial X_{3}^{K}}{\partial \widehat{X_2}} & \frac{\partial X_{3}^{K}}{\partial \widehat{X_3}}\\
  \end{array}
\right)
=\sum_{r=1}^{n_{\rm sh}} \left(
  \begin{array}{ccc}
X_{r1}^{K}\frac{\partial \widehat{\theta}_r(\widehat{\xi}_l^ F)}{\partial \widehat{X_1}} & X_{r1}^{K}\frac{\partial \widehat{\theta}_r(\widehat{\xi}_l^ F)}{\partial \widehat{X_2}} & X_{r1}^{K}\frac{\partial \widehat{\theta}_r(\widehat{\xi}_l^ F)}{\partial \widehat{X_3}} \\
X_{r2}^{K}\frac{\partial \widehat{\theta}_r(\widehat{\xi}_l^ F)}{\partial \widehat{X_1}} & X_{r2}^{K}\frac{\partial \widehat{\theta}_r(\widehat{\xi}_l^ F)}{\partial \widehat{X_2}} & X_{r2}^{K}\frac{\partial \widehat{\theta}_r(\widehat{\xi}_l^ F)}{\partial \widehat{X_3}}\\
X_{r3}^{K}\frac{\partial \widehat{\theta}_r(\widehat{\xi}_l^ F)}{\partial \widehat{X_1}} & X_{r3}^{K}\frac{\partial \widehat{\theta}_r(\widehat{\xi}_l^ F)}{\partial \widehat{X_2}} & X_{r3}^{K}\frac{\partial \widehat{\theta}_r(\widehat{\xi}_l^ F)}{\partial \widehat{X_3}}\\
  \end{array}
\right).
\end{equation}
Finally, we get 
\begin{equation}
\label{eqa:partialDT}
\frac{\partial \bigl(\widehat{\nabla} T_{K}^T(\widehat{\xi}_l^ F)_{sk}\bigr)}{\partial X_{ji}}=\sum_{r=1}^{n_{\rm sh}}\frac{\partial X_{rk}^{K}}{\partial X_{ji}}\frac{\partial \widehat{\theta}_r(\widehat{\xi}_l^F)}{\partial\widehat{X_s}}=\begin{cases}\delta_{ik}
\frac{\partial \widehat{\theta}_{\widehat{j}^{-1}_{K}(j)}(\widehat{\xi}_l^ F)}{\partial\widehat{X_s}},&\mbox{ if } j\in \widehat{j}_{K}(\{1,\ldots,n_{\rm sh}\})\\0&\mbox{ otherwise}\end{cases}.
\end{equation}
This finishes the computation of the second term. 

\subsection{\label{app:dB_dX}Computing $\frac{\partial B}{\partial X}$ }
As in \eqref{eqa:discretePDE_A} we have $B(X)_{(q,r),(k,s)}=B(e_r\theta_q,e_s\theta_k)=B_1(e_r\theta_q,e_s\theta_k)+B_2(e_r\theta_q,e_s\theta_k)$, then $\frac{\partial B(X)_{(q,r),(k,s)}}{\partial X}=\frac{\partial B_1(e_r\theta_q,e_s\theta_k)}{\partial X_{ji}}+\frac{\partial B_2(e_r\theta_q,e_s\theta_k)}{\partial X_{ji}}$.
For the first partial derivative,  one obtains
\begin{align}
\label{eqa:dB_1_dX}
\begin{split}
\frac{\partial B_1(e_r\theta_q,e_s\theta_k)}{\partial X_{ji}}=\lambda\sum_{K\in \widehat{j}^{-1}(j)_1\cap\widehat{j}^{-1}(q)_1\cap\widehat{j}^{-1}(k)_1}\sum_{l=1}^{l_q}&\Bigg\{\frac{\partial}{\partial X_{ji}}\bigl(\omega_{lK}\bigr)\nabla_r \theta_q(\xi_l^{K})\nabla_s \theta_k(\xi_l^{K})\\
&+\omega_{lK}\frac{\partial}{\partial X_{ji}}\bigl(\nabla_r \theta_q(\xi_l^{K})\bigr)\nabla_s \theta_k(\xi_l^{K})\\
                                    &+\omega_{lK}\nabla_r \theta_q(\xi_l^{K})\frac{\partial}{\partial X_{ji}}\bigl(\nabla_s \theta_k(\xi_l^{K})\bigr)\Bigg\}.
\end{split}
\end{align}
Here we use the convention that $\widehat{j}^{-1}(j)_1$ is the projection of the set $\widehat{j}^{-1}(j)$ to the first component (the index of the element).  

We thus have to compute the three partial derivatives $\frac{\partial}{\partial X_{ji}}\bigl(\omega_{lK}\bigr)$, $\frac{\partial}{\partial X_{ji}}\bigl(\nabla_r\theta_q(\xi_l^K)\bigr)$ and $\frac{\partial}{\partial X_{ji}}\bigl(\nabla_s\theta_k(\xi_l^K)\bigr)$.
We start with the first partial derivative and we obtain
\begin{equation}
\label{eqa:domega_dX}
\frac{\partial}{\partial X_{ji}}\bigl(\omega_{lK}\bigr)=\frac{\partial}{\partial X_{ji}}\bigl(\widehat{\omega_{l}}\det (\widehat{\nabla} T_K(\widehat{\xi_l})\bigr)=\widehat{\omega_l}\frac{\partial}{\partial X_{ji}}\bigl(\det( \widehat{\nabla}T_K(\widehat{\xi_l})\bigr)
\end{equation}
For an invertible matrix $A$, we have $\frac{d}{d\alpha}\det(A)=\det(A)\tr(A^{-1}\frac{dA}{d\alpha})$, so we get
\begin{equation}
\frac{\partial}{\partial X_{ji}}\bigl(\omega_{lK}\bigr)=\widehat{\omega_l}\det\left( \widehat{\nabla} T_K(\widehat{\xi_l})\right)\tr\left((\widehat{\nabla}T_K(\widehat{\xi}_l))^{-1}\frac{\partial \widehat{\nabla}T_K(\widehat{\xi}_l)}{\partial X_{ji}}\right).
\end{equation}
The partial $X_{ji}$ derivative of the matrix $\widehat{\nabla}T_K(\widehat{\xi}_l)$ has been calculated in \eqref{eqa:partialDT}, where we have to replace $\widehat{\xi}_l^F$ with $\widehat{\xi}_l$.

The second partial derivative in $\frac{\partial B_1}{\partial X_{ji}}$ is easily calculated   
\begin{align}
\label{eqa:dNtheta_dX}
\begin{split}
\frac{\partial \nabla_r \theta_q}{\partial X_{ji}}(\xi_l^K)&= \frac{\partial}{\partial X_{ji}}\left[\left(\widehat{\nabla}T_K(\hat\xi_l)^T\right)^{-1}\widehat{\nabla} \widehat{\theta}_{\widehat{j}_{K}^{-1}(q)}(\hat\xi_l)\right]_r.
\end{split}
\end{align}
We now refer to equations \eqref{eqa:dNT_dX} to \eqref{eqa:partialDT} in in Section \ref{app:dJ_dX} to conclude the computation. Here, of course, the surface quadrature point $\hat\xi^F_l$ has to be replaced by the volume quadrature point $\hat\xi_l$. The third partial derivative in \eqref{eqa:dB_1_dX} is completely analogous to the second.

For the partial derivative $\frac{\partial B_2}{\partial X}$ we obtain
\begin{align}
\label{eqa:dB_2_dX}
\begin{split}
\frac{\partial B_2(e_r\theta_q,e_s\theta_k)}{\partial X_{ji}}&=2\mu\sum_{K\in \widehat{j}^{-1}(j)_1\cap\widehat{j}^{-1}(q)_1\cap\widehat{j}^{-1}(k)_1}\sum_{l=1}^{l_q}\\
&\Bigg\{\frac{\partial}{\partial X_{ji}}\bigl(\omega_{lK}\bigr)\varepsilon\left(e_r \theta_q(\xi_l^{K}))\right):\varepsilon\left(e_s \theta_k(\xi_l^{K})\right)\\
&+\omega_{lK}\frac{\partial}{\partial X_{ji}}\bigl(\varepsilon\left(e_r \theta_q(\xi_l^{K}))\right)\bigr):\varepsilon\left(e_s \theta_k(\xi_l^{K})\right)\\
                                    &+\omega_{lK}\varepsilon\left(e_r \theta_q(\xi_l^{K}))\right):\frac{\partial}{\partial X_{ji}}\bigl(\varepsilon\left(e_s \theta_k(\xi_l^{K})\right)\bigr)\Bigg\}.
\end{split}
\end{align}
The first term is calculated in \eqref{eqa:domega_dX}. For the second term, we observe that the linear elastic strain tensor field is given by 
\begin{equation}
\varepsilon(e_r\theta_q(\xi_l^ K))_{ab}= \frac{1}{2}\left( \delta_{rb}\nabla_a \theta_q(\xi_l^ K))+\delta_{ra}\nabla_b \theta_q(\xi_l^ K)\right),
\end{equation}
and we  refer to the argument following Eq.\ \eqref{eqa:dNtheta_dX} to conclude the computation of \eqref{eqa:dB_2_dX}.

\subsection{\label{app:dF_dX} Computing $\frac{\partial F}{\partial X}$}

We have $\frac{\partial F_{(q,r)}}{\partial X_{ji}}=\frac{\partial F_{(q,r)}^{\rm vol}}{\partial X_{ji}}+\frac{\partial F_{(q,r)}^{\rm surf}}{\partial X_{ji}}$. Starting with the volume term, we obtain
\begin{align}
\label{eqa:dVolF_dX}
\begin{split}
\frac{\partial F_{(q,r)}^{\rm vol}}{\partial X_{ji}}&=\sum_{K\in \widehat{j}^{-1}(j)_1\cap\widehat{j}^{-1}(q)_1 }\sum_{l=1}^{l_q}\Bigg\{\frac{\partial}{\partial X_{ji}}\bigl(\omega_{lK}\bigr)f_r(\xi_l^K)\theta_q(\xi_l^K)\\
&+\omega_{lK}\frac{\partial}{\partial X_{ji}}\bigl(f_r(\xi_l^K)\bigr)\theta_q(\xi_l^K)
+\omega_{lK}f_r(\xi_l^K)\frac{\partial}{\partial X_{ji}}\bigl(\theta_q(\xi_l^K)\bigr)\bigg\}
\end{split}
\end{align}
The partial derivative of the volume weight has been calculated in \eqref{eqa:domega_dX}. The third term in \eqref{eqa:dVolF_dX} vanishes, as $\theta_q(\theta_l)=\hat\theta_{\widehat{j}^{-1}_K(q)}(\hat\xi_l)$ does not depend on $X_{ji}$. Unless the volume force density $f$ does depend explicitly on the position (like in the case of centrifugal loads) this term vanishes as the third one. 

Finally we have to calculate the partial derivative of the surface loads
\begin{align}
\label{eqa:dSurF_dX}
\begin{split}
\frac{\partial F_{(q,r)}^{\rm surf}}{\partial X_{ji}}&=\sum_{F\in \mathcal{N}_h\atop K(F)\in \widehat{j}^{-1}(j)_1\cap\widehat{j}^{-1}(q)_1}\sum_{l=1}^{l_q^F}\Bigg\{\frac{\partial}{\partial X_{ji}}\bigl(\omega_{lF}\bigr)f_r(\xi_l^F)\theta_q(\xi_l^F)\\
&+\omega_{lF}\frac{\partial}{\partial X_{ji}}\bigl(f_r(\xi_l^F)\bigr)\theta_q(\xi_l^F)
+\omega_{lF}f_r(\xi_l^F)\frac{\partial}{\partial X_{ji}}\bigl(\theta_q(\xi_l^F)\bigr)\Bigg\}.
\end{split}
\end{align}
The first term can be calculated with the aid of \eqref{eqa:domega_F_dX}. For the second and third term, the same reasoning applies as for partial derivative of the volume force.


\small


   \begin{thebibliography}{marke}
    \bibitem{AGN} S. Agmon, A. Douglis and L. Nirenberg, \textit{Estimates Near the Boundary for Solutions of
    Elliptic Partial Differential Equations Satisfying General Boundary Conditions II}, Communications On Pure And Applied Mathematics, Vol. XVII,
    pp. 35-92, 1964.
    
    \bibitem{All} G. Allaire, \textit{Numerical Analysis and Optimization}, Oxford University Press, Oxford (2007)
    \bibitem{ABFJ} G.\ Allaire, E.\ Bonneter, G.\ Francfort and F.\ Jouve,  \textit{Shape Optimization by the Homogenization Method}, Numer. Math., {\bf 76}, 1997, 27-68.
    
    \bibitem{BSSST} I. Babu$\check{\rm s}$ka, Z. Sawlan, M. Scavino B. Szab\'o and R. Tempone, Spatial Poisson processes for crack initiation, (2018) arXiv:1805.03433.

   \bibitem{BR} W. Bangerth, R. Rannacher, \textit{Adaptive Finite Element Methods for Differential Equations}, Birkh\"auser 2003. 
   
   \bibitem{BXZ} R. E. Bank, J.-C. Xu and B. Zheng, Superconvergent derivative recovery 	for Lagrange triangular elements of degree p on unstructured grids, SIAM J. Numer.\ Anal. Vol. 45,  2007, 2032-2046.

    \bibitem{BHR} M. B\"aker, H. Harders and J. R\"osler, \textit{Mechanical Behaviour of Engineering Materials: Metals, Ceramics,
Polymers, and Composites}, German edition published by Teubner
Verlag (Wiesbaden, 2006), Springer, Berlin Heidelberg 2007.

\bibitem{BC} S. B. Batdorf and J. G. Crosse,\textit{ A statistical theory for the fracture of brittle structures subject to nonuniform
polyaxial stress}, J. Appl. Mech. 41, 459--465 (1974)

\bibitem{BHGS} M. Bolten, C. Hahn, H. Gottschalk and M. Saadi, Numerical shape optimization to decrease the failure probability of ceramic structures, (2017) arXiv:1705.05776. 

\bibitem{BGS} M. Bolten, H. Gottschalk and S. Schmitz, \textit{Minimal failure probability for ceramic design via shape control}, J. Optim. Theory  Appl.,
{\bf 166} (2015), 983--1001.

	\bibitem{BACF} L.P. Borrego, L.M. Abreu, J.M. Costa and J.M. Ferreira, \textit{Analysis of low cycle fatigue in AlMgSi aluminium alloys}, Engineering Failure Analysis {\bf 11} (2004) 715-725.

    \bibitem{BS} A.\ Borzi and V.\ Schulz, \textit{Computational Optimization of Systems governed by Partial Differential Equations}, SIAM series on computational engineering, SIAM 2012.
    
    \bibitem{Bra} D. Braess, \textit{Finite Elemente - Theorie, schnelle L\"oser und Anwendungen in der Elastizit\"atstheorie}, fourth edition, Springer, Berlin, 2007.

    \bibitem{Cia} P. Ciarlet, \textit{Mathematical Elasticity - Volume I: Three-Dimensional Elasticity}, Studies in Mathematics and its Applications, Vol. 20, North-Holland, Amsterdam, 1988.
    
    \bibitem{EG} A. Ern and J.-L. Guermond, \textit{Theory and Practice of Finite Elements}, Springer, New York, 2004.
 
     \bibitem{EM} L. A. Escobar and W. Q. Meeker, \textit{Statistical Methods for Reliability Data}, Wiley-Interscience Publication, New York, 1998.
     
     \bibitem{FM} D. Munz and T. Fett, \textit{Ceramics}, (engl. edition) Springer Verlag Berlin, 2012.
     \bibitem{GSDKS} H. Gottschalk, M. Saadi, O. T. Doganay, K. Klamroth and S. Schmitz, Adjoint Method to Calculate Shape Gradients of Failure Probabilaties for Turbomachinery Components, ASME-Turbo-Expo GT2018-75759 (2018). 
     \bibitem{GS} H. Gottschalk and S. Schmitz, Optimal Reliability in Design for Fatigue Life I: Existence of optimal shapes, SIAM J. Control Optim. {\bf 52} Vol. 5, (2014), 2727-2752.
     
     
     \bibitem{GSSKRB} H. Gottschalk,  S. Schmitz, T. Seibel, R. Krause,
G. Rollmann and T. Beck, \textit{Probabilistic Schmid Factors and Scatter of LCF Life} ,  Materials Science and Engineering 46 (2) (2015), 156-164..

     \bibitem{Got} G. Gottstein, \textit{Physical foundations of Material Science}, Springer 2004.
     
    \bibitem{HM} J. Haslinger and R. A. E. M\"akinen, \textit{Introduction to Shape Optimization - Theory, Approximation and Computation}, SIAM - Advances in Design and Control, 2003.
    
    \bibitem{HV} Hertel, O., Vormwald, M., \textit{Statistical and geometrical size effects in notched members based on
weakest-link and short-crack modelling}, Engineering Fracture Mechanics, \textbf{95} (2012), pp 72--83.
    
      \bibitem{HS} M. Hoffmann and T. Seeger, \textit{A Generalized Method for Estimating Elastic-Plastic Notch Stresses and Strains, Part 1: Theory}, Journal of Engineering Materials and Technology, 107, pp. 250-254, 1985.
      
        \bibitem{Kal} O. Kallenberg, \textit{Random Measures}, Akademie-Verlag, Berlin 1983.
        
    \bibitem{KJW} M. Knop, R. Jones, L. Molent and L. Wang,
            \textit{On Glinka and Neuber methods for calculating notch tip strains under cyclic load spectra},
            International Journal of Fatigue, Vol. 22, (2000) 743--755.
            
    \bibitem{LNT} W. B. Liu, P. Neittaanm\"aki and D. Tiba, \textit{Existence for shape optimization problems in arbitrary dimension}, SIAM J.\ Contr.\ Optimization {\bf 41} (2003) 1440-1454.
    
    \bibitem{LLMMW} G. R. Leverant, D. L. Littlefield,
R. C. McClung, H. R. Millwater, and J. Y. Wu, \textit{A Probabilistic Approach to Aircraft Turbine Rotor Material Design},  Paper 97-GT-22, ASME Turbo Expo '97, Orlando, Florida,
June 1997.
\bibitem{LSGB} L. M\"ade, S. Schmitz, H. Gottschalk and T. Beck, Combined notch and size effect modeling in a local probabilistic
approach for LCF, Comp.\ Materials Science 142 (2018) 377–388.
    \bibitem{Neu} H. Neuber, \textit{Theory of Stress Concentration for Shear-Strained Prismatical Bodies with Arbitrary Nonlinear Stress-Strain Law}, J. Appl. Mech. 26, 544, 1961.
    
     \bibitem{RV} D. Radaj and M. Vormwald, \textit{Erm\"udungsfestigkeit}, third edition, Springer, Berlin Heidelberg, 2007.
     
    \bibitem{RO} W. Ramberg and W. R. Osgood, \textit{Description of Stress-Strain Curves by Three Parameters}, Technical Notes - National Advisory Committee For Aeronautics, No. 902, Washington DC., 1943.
    \bibitem{RBZ} H. Riesch-Oppermann, A. Br\"uckner-Foit, C. Ziegler, \textit{ STAU - a general purpose tool for probabilistic
reliability assessment of ceramic components under multi axial loading}, in: Proceedings of the 13th
International Conference on ECF 13, San Sebastian (2000)
    
    \bibitem{Schmi} S. Schmitz, \textit{A Probabilistic Local Model for Low Cycle Fatigue -- New Aspects of Structurel Mechanics}, Dissertation Lugano and Wuppertal 2014, appeared in Hartung-Gorre Verlag, 2014.
    
    \bibitem{SRGK1}  S. Schmitz, G. Rollmann, H. Gottschalk  and R.
Krause, \textit{Risk estimation for LCF crack initiation}, Proc.
ASME Turbo Expo 2013, GT2013-94899, arXiv:1302.2909v1.

    \bibitem{SRGK2} S. Schmitz,  G. Rollmann, H. Gottschalk
and R.\ Krause, \textit{Probabilistic analysis of the LCF crack initiation life for a turbine blade under
thermo-mechanical loading}, Proc.\ Int.\ Conf.\ LCF 7 (September 13).

        \bibitem{SSBRKG} S. Schmitz, T. Seibel, T. Beck, G. Rollmann, R. Krause and H. Gottschalk, \textit{A probabilistic Model for LCF}, Computational Materials Science {\bf 79} (2013), 584-590.
        
        \bibitem{SZ} J. Sokolowski and J.-P. Zolesio, \textit{Introduction to Shape Optimization -  Shape Sensivity Analysis}, first edition, Springer, Berlin Heidelberg, 1992.
        
        \bibitem{SMB} D. Sornette, T. Magnin, and Y. Brechet, \textit{The Physical Origin of the Coffin-Manson Law in
Low-Cycle Fatigue}, Europhys. Lett., 20 (1992), pp. 433--438.

\bibitem{Troe} F. Tr\"oltzsch, Optimal Control of Partial Differential Equations (in German), Vieweg + Teubner, 2010.

 \bibitem{Wat} S. Watanabe, \textit{On discontinuous additive functionals and Lévy measures of a Markov process}, Japan. J. Math. 34 (1964).
 
 \bibitem{Wei} E. W. Weibull, \textit{A Statistical Theory of the Strength of Materials}, Ingeniors Vetenskaps Akad. Handl.
151, 1--45 (1939)
 
 \bibitem{Wit} K. Wittig, \textit{Construction of a gas turbine for model air planes (in German)}, Munich 1993, www.calculix.de.
  \end{thebibliography}
\end{document}